\documentclass{article}

\usepackage[margin= 1in]{geometry}
\usepackage{amsmath,amsthm,amssymb,bm,graphicx }
\usepackage{subfig}
\usepackage{fancyvrb}
\usepackage{mathtools}
\usepackage{multirow}
\usepackage[hyphens]{url}
\usepackage{hyperref,xcolor, soul}
\sethlcolor{pink}
\usepackage{tikz}
\usepackage{comment}
\usetikzlibrary{arrows,matrix,positioning}
\usepackage{enumitem}
\usepackage{verbatim}
\usepackage{accents}
\usepackage{lipsum}
\usepackage{adjustbox}
\usepackage[colorinlistoftodos]{todonotes}
\usepackage{tabularray}
\usepackage{tikz-cd}
\usepackage{authblk}

\newcommand{\F}{\mathbb{F}}
\newcommand{\Q}{\mathbb{Q}}
\newcommand{\Z}{\mathbb{Z}}
\newcommand{\N}{\mathbb{N}}

\newcommand{\p}{\mathfrak{p}}
\newcommand{\OO}{\mathcal{O}}
\newcommand{\inv}{^{-1}}
\DeclareMathOperator{\Gal}{Gal}

\newtheorem{theorem}{Theorem}[section]
\newtheorem{lemma}[theorem]{Lemma}
\newtheorem{proposition}[theorem]{Proposition}

\newtheorem{corollary}[theorem]{Corollary}

\newtheorem*{mainthm}{Theorem}
\theoremstyle{definition}
\newtheorem{definition}{Definition}[section]

\theoremstyle{remark}
\newtheorem*{remark}{Remark}

\hypersetup{
	colorlinks=true,
	linkcolor=blue,
	citecolor=blue,
	linktoc=page}

\title{On Classifying Extensions of $p$-adic Fields}
\date{\vspace{-5ex}}
\author{Shreya Dhar\thanks{University of Chicago, \texttt{shreyadhar@uchicago.edu}}, 
River Newman\thanks{Yale University, \texttt{river.newman@yale.edu}}, 
Grayson Plumpton\thanks{University of Toronto, \texttt{grayson.plumpton@mail.utoronto.ca}}, 
Chenglu Wang\thanks{University of Pennsylvania, \texttt{chengluw@sas.upenn.edu}}}

\begin{document} 

\maketitle
\begin{abstract}
    Let $p$ be a prime and let $\Q_p$ be the field of $p$-adic numbers. It is known that the finite extensions of $\Q_p$ of a given degree are finite up to isomorphism. Given a cubic field extension $L$ of $\Q_p$ generated by the root of an irreducible polynomial $h$, we present a practical (closed-form) method to determine the isomorphism class in which $L$ lives, based on the coefficients of $h$. We discuss the subtleties of the wildly ramified case, when the degree of the extension coincides with $p$, the characteristic of the residue field. We also present a method for tamely ramified extensions of arbitrary prime degree.
\end{abstract}

\section{Introduction}
\subsection{Notation}

Throughout we will use the following notation
\begin{center}
	\begin{tabular}{cll}
        $p$ & : & a prime number\\
		$K$ & : & an unramified extension of the $p$-adic numbers $\Q_p$, with finite degree; \\
		$v_p$ & : & the $p$-adic valuation on $K$ (so $v_p(p) = 1$); \\
		$|\cdot|_p$ & : & the norm induced by $v_p$, $|x|_p = p^{-v(x)}$; \\
		$\OO_K$ & : & the valuation ring or ring of integers of $K$; \\
        $C_r$ & : & the cyclic group of order $r$.
	\end{tabular}
\end{center}

We may drop the subscripts on notation when $p$ is obvious or assumed from context.

\subsection{Background and overview}

The study of field extensions of $p$-adic fields is of interest in number theory and crucial to  understanding classical groups over $p$-adic fields. Morris in \cite{Morris} showed that, for such groups of degree $n$ over $\Q_p$, the conjugacy classes of maximal elliptic tori are determined by isomorphism classes of extensions of $\Q_p$ with degrees up to $n$. Chinner in \cite{Chinner} uses this framework to classify 4-dimensional special orthogonal groups over $p$-adic fields. This connection makes the classification of $p$-adic extensions an essential step for understanding these tori and constructing supercuspidal representations, which play a key role in the local Langlands correspondence.

Previous results, such as \cite{Pauli} and \cite{Amano} (for tamely and wildly ramified extensions respectively; results summarized helpfully in \cite{Database}), have effectively classified all such extensions. However, given an irreducible polynomial $f(x)$ of degree $n$, we are not aware of a simple procedure by which one can in general determine which isomorphism class the extension $f(x)$ generates will belong to.

Such a procedure exists for quadratic extensions, and is a simple result following directly from the quadratic formula: given $f \in \Q_p[x]$ quadratic with discriminant $\Delta$, $\Q_p[x]/\langle f(x)\rangle \cong \Q_p(\sqrt{\Delta})$, whose isomorphism class depends only on $\Delta$'s equivalence class up to squares in $\Q_p^\times$—easily calculable outside of the wildly ramified case $p = 2$, and for $p=2$, it is only moderately more complicated. For degrees $>2$, though, the state of affairs is more difficult; even for cubic and quartic polynomials, the relevant formulae do not offer much help.
    
In this paper, we develop a similar procedure for cubic extensions, which, given an irreducible cubic polynomial over $\Q_p$, allows the isomorphism class of the corresponding extension to be determined by the cube closed-form expression depending only on the coefficients of the polynomial and $p$.

In the case of tame ramification only (i.e., when $p \nmid n$), we extend this result to obtain a similar procedure for polynomials generating totally tamely ramified extensions of arbitrary degree $n$, over a base field $K$ which can be any unramified extension of $\Q_p$. Since all extensions of $p$-adic fields decompose uniquely into totally ramified and unramified parts, this is likely applicable to extensions which are tamely, but not totally, ramified.

We hope that these results may be of use in the study of low-degree extensions of $p$-adic fields, and possibly also the structure theory of certain classical groups over $\Q_p$.

\subsection{Main results}

\begin{mainthm}[\textbf{A}] (Theorem \ref{thm:tame-prime-classification})
    Let $p,q$ distinct primes, $K/\Q_p$ unramified of degree $m$. Let $\zeta \in K$ a primitive $(p^m - 1)$th root of unity, and $d = \gcd(q, p^m-1)$.
    
    Suppose $f(x)\in \OO_K[x]$ an irreducible polynomial of degree $q$, $\theta \in \overline{K}$ a root of $f$.

    If $q(q-1) \mid v(\Delta_f)$, $K(\theta)/K$ is the unique unramified extension of degree $q$ over $K$.
    
    Otherwise, $K(\theta) \cong K(\sqrt[q]{\zeta^r p})$ for the (unique mod $d$) value of $r$ such that $\zeta^r u^\ell$
    is an $q$th power mod $p$, where $\ell$ is the inverse mod $q$ of $\frac{v(\Delta_f)}{q-1}$ and $u$ is the unit part of $\Delta_f$.
\end{mainthm}

\begin{mainthm}[\textbf{B}] (Theorems \ref{Wild Classify Thm 1}, \ref{main theorem Galois}, \ref{main theorem non-Galois})
   Let $f(x) \in \Q_3[x]$ be an irreducible cubic polynomial. Then we can classify up to isomorphism the extension generated by $f(x)$ using the following procedure:

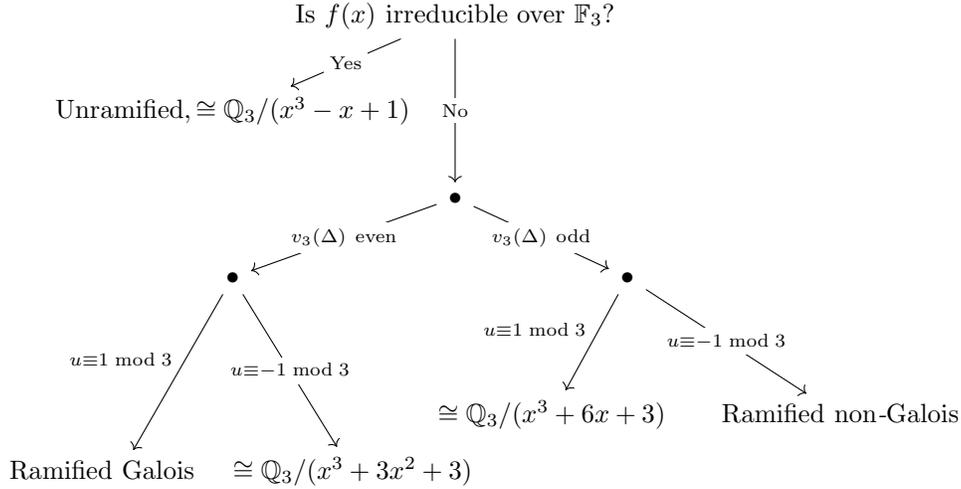
\begin{center}

  \[\begin{tikzcd}
 &[-100pt]&[-100pt]&[-100pt]\mathrm{Is\ } f(x) \mathrm{\ irreducible \ over \ } \F_3? \ar[dll,"\mathrm{Yes}"{description}]\ar[dd,"\mathrm{No}"{description}]&[-100pt]&[-50pt]\\
 &\mathrm{Unramified}, \cong \Q_3 / ( x^3-x+1) && {} \\
 &&&
 \bullet\ar[dll,"v_3(\Delta) \mathrm{ \ even}"{description}]\ar[drr,"v_3(\Delta) \mathrm{ \ odd}"{description}]\\
 &\bullet \ar[dddl,"u\equiv 1 \bmod 3",swap] \ar[dddr,"u\equiv -1 \bmod 3"{description}]&&& &\bullet\ar[ddl,"u\equiv 1 \bmod 3",swap]\ar[ddr,"u\equiv -1 \bmod 3"{description}]\\
 &&&&&\\
 &{}&&&\cong \Q_3 / (x^3+6x+3) &&\mathrm{Ramified \ non\text{-}Galois}\\[-15pt]
  \mathrm{Ramified\  Galois \ \ \ \ } && \cong \Q_3 / (x^3+3x^2+3)
\end{tikzcd}\] 
\end{center}

where $\Delta$ is the discriminant of $f(x)$ and $u \in \Z_3^\times$ is the unit part of $\Delta$.

\bigskip

In the ramified Galois case, transform $f$ into monic depressed form $f(x) = x^3 + \alpha x + \beta$. Let $u$ the unit part of $\beta$, $w$ the unit part of $\alpha$ (if $\alpha\neq 0$).

Then $\Q_3(\theta)/\Q_3 \cong \Q_3[x]/( x^3-3x^2+3\tau)$ for the unique $\tau \in \{1,4,7\}$ such that $\pm \tau \equiv t$ mod 9, where
\begin{align*}
    t = \begin{cases}
        w + \frac{1-u}{3} & v_3(\beta) \equiv 0 \bmod 3,\quad u \equiv 1 \bmod 9\\
        w + \frac{1+u}{3} & v_3(\beta) \equiv 0 \bmod 3,\quad u \equiv -1 \bmod 9\\
        u & v_3(\beta) \equiv 1 \bmod 3\\
        u^2 & v_3(\beta) \equiv 2 \bmod 3
    \end{cases}
\end{align*}
(In particular, $3\mid v_3(\beta) \implies \alpha\neq 0$.)

\bigskip

Finally, in the ramified non-Galois case, again take the monic depressed form $f(x) = x^3 + \alpha x + \beta$. Let $u$ the unit part of $\beta$, $w$ the unit part of $\alpha$ (if $\alpha\neq 0$). Let $r = v_3(\alpha)$, $m = \lfloor v_3(\beta)/3 \rfloor$. Then
\begin{enumerate}
    \item If $v_3(\beta) \equiv 0 \bmod 3$, $\Q_3(\theta)/\Q_3 \cong \Q_3[x]/(x^3 + 3x + 3)$.

    \item If $v_3(\beta) \equiv 1 \bmod 3$ and $r = 2m + 1$, $\Q_3(\theta)/\Q_3 \cong \Q_3[x]/(x^3 + 3x + 3)$.

    \item If $v_3(\beta) \equiv 1 \bmod 3$ and $r > 2m + 1$, $\Q_3(\theta)/\Q_3 \cong \Q_3[x]/(x^3 + 3\tau)$ for the unique $\tau \in \{1,4,7\}$ such that
    \[
        \pm\tau \equiv \frac{u}{1 - 3^{r - 2m - 1} w} \equiv \begin{cases}
            \frac{u}{1-3w} & r = 2m + 2\\
            u & r > 2m + 2
        \end{cases} \quad \mod 9
    \]

    \item If $v_3(\beta) \equiv 2 \bmod 3$, $\Q_3(\theta)/\Q_3 \cong \Q_3[x]/(x^3 + 3\tau)$ for the unique $\tau \in \{1,4,7\}$ such that
    \[
        \pm\tau \equiv \frac{1 + w(-3)^{r-2m-1}}{u} \equiv \begin{cases}
            \frac{1-3w}{u} & r = 2m + 2\\
            \frac{1}{u} & r > 2m + 2
        \end{cases} \quad\mod 9
    \]
\end{enumerate}

\end{mainthm}

\section{Tame ramification in degree $n$}

We give general results on tamely ramified extensions of the base field $K$, itself an unramified finite extension of $\Q_p$. In particular, in the case where $n$ prime and $\neq p$, we provide a closed-form method to determine up to isomorphism the extension generated by any irreducible degree $n$ polynomial $f \in \OO_K[x]$.

\subsection{$k$-Eisenstein polynomials}

\begin{definition}
    Consider a polynomial $h(x) = \sum_{i=0}^n a_i x^i$ in $\OO_K[x]$, $a_n \neq 0$. For some $k$ with $1 \leq k < n$, $h$ is $k$-Eisenstein
    \footnote{
    These are called generalized Eisenstein polynomials by some other sources.
    }
    if $v(a_n) = 0$, $v(a_0) = k$, and $v(a_i) + ik/n \geq k$ for all $i$.
    Equivalently, $h$ is $k$-Eisenstein if $v(a_n) = 0$ and its Newton polygon consists of a single segment of slope $-k/n$.
\end{definition}
Note that since $v(a_n) = 0$, $a_n$ a unit in $\OO_K$. So, without loss of generality, we can take $h(x)$ to be monic (without losing the fact that it is in $\OO_K[x]$).

\begin{proposition}
    Let $h \in \OO_K[x]$ with degree $n$, $p\nmid n$, and $k$-Eisenstein for some $k$ with $(k,n) = 1$. If $\theta$ is a root of $h$, $K(\theta)/K$ is degree $n$ and totally ramified. In particular, $v(\theta) = k/n$.
\end{proposition}
\begin{proof}
    By the condition on the Newton polygon, $v(\theta) = k/n$. By Euler's theorem, since $(k,n) = 1$, we have $v(\theta^{k^{\varphi(n)-1}}) = k^{\varphi(n)}/n = j + 1/n$ for some $j \in \Z$, and therefore the ramification index of $K(\theta)/K$ is divisible by $n$. As the degree is at most $n$, this means the degree and ramification index are both $n$ exactly.
\end{proof}
\begin{corollary}
    Any $k$-Eisenstein polynomial of degree $n$ such that $(k,n) = 1$ is irreducible.
\end{corollary}

\begin{proposition}\label{prop:min-poly-k-eis}
    If $\theta$ has degree $n$ over $K$ and $v(\theta) = k/n$ where $1 \leq k < n$, $(k,n) = 1$, the minimal polynomial of $\theta$ over $K$ is $k$-Eisenstein (and, as a prerequisite, is in $\OO_K[x]$).
\end{proposition}
\begin{proof}
    The minimal polynomial of $\theta$ over $K$ yields (via multiplication by a common denominator) a nontrivial $\OO_K$-linear relation $\sum_{j=0}^n a_j \theta^j = 0$, where $a_j \in \OO_K$, $a_n \neq 0$, and $\min_j v(a_j) = 0$.
    
    By the ultrametric inequality, it follows that the maximum norm of any term $\max_j |a_j \theta^j|$ occurs at least twice among all the terms. Equivalently, the minimum of $v(a_j \theta^j)$ occurs at least twice.
    
    But $v(a_j \theta^j) = jk/n + v(a_j)$, where $v(a_j)$ an integer. Since $(k,n) = 1$, these all have distinct fractional part and cannot be equal, except for $j=0,n$. Therefore $v(a_0) = v(a_n \theta^n) = v(a_n) + k$. Furthermore, since this valuation must be the minimum among all terms, $ik/n + v(a_i) = v(a_i \theta^i) \geq v(a_n) + k$ for all $i$. In particular, this means $v(a_i) \geq v(a_n)$. Since by assumption $\min_j v(a_j) = 0$, $v(a_n) = 0$. Then $ik/n + v(a_i) \geq k$ for all $i$.
    
    Therefore the polynomial $\sum_{j=0}^n a_j x^j \in \OO_K[x]$ is $k$-Eisenstein and has root $\theta$. Let $h(x) = a_n\inv \sum_{j=0}^n a_j x^j$; since $v(a_n) = 0$, this is in $\OO_K[x]$ and is also $k$-Eisenstein, thus irreducible. Since $h(x)$ is monic, it is the minimal polynomial of $\theta$.
\end{proof}

\begin{proposition}\label{prop:k-eis-shift}
    Let $K(\theta)/K$ a totally ramified extension of degree $n$. There exists $t \in K$, $r \in \Z$, $1 \leq k < n$ such that $(k,n) = 1$ and the minimal polynomial of $p^r(\theta + t)$ is $k$-Eisenstein. Furthermore, $k$ is unique.
\end{proposition}
\begin{proof}
    Let $\pi$ be any uniformizer of $K(\theta)$. As $K(\theta)/K$ totally ramified, $v(\pi) = 1/n$ and $K(\theta) = K(\pi)$. Then we can represent $\theta$ in terms of a power basis, $\theta = \sum_{j=0}^{n-1} c_j \pi^j$  where $c_j \in K$.
    
    Take $t = -c_0$, $r = -\min_i v(c_i) = -\lfloor v(\theta) \rfloor$. Then
    \[
    p^r(\theta + t) = \sum_{j=1}^{n-1} p^r c_j \pi^j
    \]
    where $v(p^r c_j \pi^j) = j/n +  v(c_j) + r$. In particular, $v(c_j) + r \geq 0$, so every term $p^r c_j \pi^j$ has positive valuation with fractional part $j/n$. Furthermore, for all $j$ minimizing $v(c_j)$, we have  $v(p^r c_j \pi^j) = j/n$.
    
    Since their fractional parts $j/n$ all differ, $v(p^r c_j \pi^j)$ are pairwise unequal. Thus by the ultrametric inequality
    \[
    v(p^r(\theta + t)) = \min_j v(p^r c_j \pi^j) = k/n
    \]
    for the least $k$ such that $v(c_k) = r$. So the ramification index of $K(p^r(\theta + t))$ is $n/(k,n)$. Since $K(p^r(\theta + t)) = K(\theta)$ is totally ramified of degree $n$, it follows that $(k,n) = 1$, and then by Proposition \ref{prop:min-poly-k-eis} the minimal polynomial of $p^r(\theta + t)$ is $k$-Eisenstein.

    For the uniqueness of $k$, suppose that $p^r(\theta + t)$, $p^s(\theta + u)$ are respectively $k$- and $j$-Eisenstein. Then $v(\theta + t) = -r + k/n$, $v(\theta + u) = -s + j/n$. If $j\neq k$, these cannot be equal, and thus
    \[
    v(t - u) = v(\theta + t - (\theta + u)) = \min\left( - r + \frac{k}{n}, - s + \frac{j}{n} \right)
    \]
    a noninteger, which contradicts the fact that $t-u \in K$. So $j = k$ and $k$ unique.
\end{proof}
\begin{corollary}
    Under the hypotheses of the above, let $f(x) = \sum_{j=0}^n a_j x^j \in \OO_K[x]$ an irreducible polynomial of which $\theta$ is a root.
    
    \begin{enumerate}
        \item The value of $t$ obtained in the above proof is the unique $t$ such that $f(x-t)$ is a depressed polynomial (lacks a term of degree $n-1$), i.e. $t = \frac{a_{n-1}}{n a_n}$.

        \item The value of $r$ is uniquely determined by this choice of $t$.
    \end{enumerate}

\end{corollary}
\begin{proof}
    \begin{enumerate}
        \item Represent $\theta + t = \sum_{j=1}^{n-1} c_j \pi^j$ as a linear transformation with respect to the power basis. Its matrix has zero diagonal, since there is no $\pi^0$ term. Examining the cofactor decomposition of the determinant shows that the characteristic polynomial has no $x^{n-1}$ term, i.e. it is depressed.

        This characteristic polynomial is a monic degree $n$ polynomial with root $\theta + t$, making it necessarily the minimal polynomial of of $\theta + t$, and thus equal to a constant multiple of $f(x-t)$. So $f(x-t)$ is depressed.

        \item Let $d_0, \ldots, d_n$ be the coefficients of $g(x)$, the minimal polynomial of $p^r(\theta + t)$. Then $v(d_0/d_n) = v(d_0) = k$ since it is monic and $k$-Eisenstein.

        $f(xp^{-r} - t) = f(p^{-r}(x - p^r t))$ is a polynomial of degree $n$ which has as a root $p^r(\theta + t)$; therefore it is equal to a scalar multiple of $g$. In particular, the ratio of its $0$th and $n$th coefficients will be equal to $d_0/d_n$. Computing this ratio and taking its valuation shows that
        \[
        k = v(d_0/d_n) = v(a_0') + nr - v(a_n)
        \]
        \[
        r = \frac{k + v(a_n) - v(a_0')}{n}
        \]
        where
        \[
        a_0' = a_0 + \sum_{j=1}^n a_j (-t)^j = a_0 + \sum_{j=1}^n a_j \left( -\frac{a_{n-1}}{n a_n} \right)^j
        \]
    \end{enumerate}
\end{proof}
\begin{remark}
   Other possible values of $t$ may exist. For example, the polynomial $f(x) = x^2 + 5x + 5 \in \Z_5[x]$ is 1-Eisenstein, and so is the shift $f(x - 5/2) = x^2 + 45/4$; the above gives the latter. But we will only consider the unique depressed case, since it is easily computable from the form of the original polynomial.
\end{remark}

\subsection{Classification of tamely totally ramified extensions}

For the following, we assume that $p\nmid n$.

\begin{theorem}
    Let $h(x) = \sum_{i=0}^n a_i x^i \in \OO_K[x]$. Suppose that $h$ is $k$-Eisenstein for some $k$ with $(k,n) = 1$. If $\theta \in \overline{K}$ is a root of $h$, then there exists $\gamma \in K(\theta)$ such that $\gamma^n = {-a_0/a_n}$ and $K(\theta) = K(\gamma)$.
    \footnote{This result and its proof are adapted and generalized from the one given in \cite{mathstack}.}
\end{theorem}
\begin{proof}
    Let $\pi \in \overline{K}$ be some fixed $n$th root of $p$. Let $v$ denote the $p$-adic valuation on $K(\pi,\theta)$ with $v(p) = 1$, i.e. extending the valuation on $K$. Then
    \[
        h(\pi^k x) = \sum_{i=0}^n a_i \pi^{ki} x^i = p^k \sum_{i=0}^n b_i x^i
    \]
    where $b_i = a_i \pi^{k(i-n)}$. For all $i$, $v(b_i) = v(a_i \pi^{ki}) - k = v(a_i) + ik/n - k \geq 0$ since $h$ is $k$-Eisenstein. If $0 < i < n$, we have also that $n\nmid ik$, since $(k,n)=1$, so $v(b_i)$ is not an integer and therefore $v(b_i)\geq 1/n$. Furthermore, $v(b_n) = v(a_n) + nk/n - k = 0$ and $v(b_0) = v(a_0) -k = 0$. 
    
    $\pi^{-k} \theta$ has valuation 0 and is a root of $h(\pi^k x)$, thus also of the polynomial
    \[
    g(x) = p^{-k}h(\pi^k x) = \sum_{i=0}^n b_i x^i
    \]
    Since each $b_i$ with $0 < i < n$ has $v(b_i) \geq 1/n$, $b_i$ is divisible by $\pi$ in $\OO_{\overline{K}}$. So $g(x) \equiv b_n x^n + b_0 \bmod \pi$ in $\OO_{\overline{K}}[x]$.

    Since $\pi^{-k}\theta$ solves $g(x)$, it solves it mod $\pi$, and furthermore $nb_n(\pi^{-k}\theta)^{n-1} \not\equiv 0 \bmod \pi$, since $v(\pi^{-k}\theta) = 0$ and $p\nmid n$. So by Hensel's lemma there exists some $\gamma' \in \overline{K}$ a root of $b_n x^n + b_0$ such that $\pi^{-k} \theta \equiv \gamma' \bmod \pi$ in $\OO_{\overline{K}}$. Note that $v(b_n) = v(b_0) = 0$, so $v(\gamma') = 0$.

    Furthermore, the fact that $\pi^{-k} \theta \equiv \gamma' \bmod \pi$ implies $v(\pi^{-k} \theta - \gamma') > 0$ and thus $v(\theta - \pi^k\gamma') > k/n$.

    The minimal polynomial of $\pi^k\gamma'$ over $K$ is $b_n x^n + p^k b_0$. Therefore all its Galois conjugates are of the form $\pi^k\gamma'\zeta^r$ where $\zeta \in \overline{K}$ is a primitive $n$th root of unity and $0 \leq r < n$.

    For any $0 < r < n$, $1 - \zeta^r$ is an integer in the unramified extension $K(\zeta)$, so $v(1 - \zeta^r)$ is a nonnegative integer. And since the roots of unity of $K(\zeta)$ correspond directly with those of its residue field, $\overline{\zeta}$ is also a primitive $n$th root of unity in $\OO_{K(\zeta)}/p\OO_{K(\zeta)}$, so $1 - \zeta^r \not\equiv 0 \bmod p$, and thus $v(1 - \zeta^r) < 1$, $v(1 - \zeta^r) = 0$. Now
    \[
        v(\pi^k\gamma' - \pi^k\gamma'\zeta^r) = k/n + v(1 - \zeta^{r}) = k/n
    \]
    So $v(\theta - \pi^k\gamma') > v(\pi^k\gamma' - \pi^k\gamma'\zeta^r)$. Therefore by Krasner's lemma $K(\pi^k\gamma') \subseteq K(\theta)$.
    
    $v(\pi^k\gamma') = k/n$, so $K(\pi^k\gamma')/K$ has ramification index $n$, thus degree $\geq n$; it follows that $K(\pi^k\gamma') = K(\theta)$. And by construction $\gamma = \pi^k\gamma'$ is an $n$th root of $-p^k b_0/b_n = -a_0/a_n$, giving the desired result.
\end{proof}
\begin{corollary}\label{cor:k-eis-nth-root-form}
    For any $\ell \in \Z$ such that $k\ell \equiv 1 \bmod n$ (since $(k,n) = 1$, such $\ell$ exists by Euler's theorem), 
    \[
    K(\theta) \cong K\left(\sqrt[n]{p(-w/a_n)^\ell}\right)
    \]
    where $w$ is the unit part of $a_0$.
\end{corollary}
\begin{proof}
    Take $\gamma$ provided by the above theorem. Then let $\phi = p^{\frac{1-k\ell}{n}} \gamma^\ell$.
    
    Clearly, $K(\phi) \subseteq K(\gamma) = K(\theta)$, and
    \[    
    \phi^n = p^{1-k\ell} \gamma^{n\ell} =  p^{1-k\ell}  (-a_0/a_n)^\ell = p(-p^{-k} a_0/a_n)^\ell = p(-w/a_n)^\ell
    \]
    so $\phi$ is an $n$th root of $p(-w/a_n)^\ell$. It follows that it has valuation $1/n$, so the ramification index of $K(\phi)/K$ is $n$, and therefore $K(\phi) = K(\theta)$. 
\end{proof}

\begin{proposition}\label{prop:nth-root-extension-isomorphism-criterion}
    Let $a, b \in \OO_K$ nonzero with $v(a) = v(b)$ and $K(\sqrt[n]{a})$, $K(\sqrt[n]{b})$ each totally ramified of degree $n$. Then $K(\sqrt[n]{a}) \cong K(\sqrt[n]{b})$ iff $a/b \in \OO_K$ is an $n$th power mod $p$.
\end{proposition}
\begin{proof}
    Suppose $\overline{a/b} \equiv \bar{c}^n \bmod p$. Then by Hensel (since $p\nmid n$) we can pick a lift $c \in \OO_K$ of $\bar{c}$ such that $a/b = c^n$. Then $\left(\sqrt[n]{a}/c\right)^n = a/c^n = b$. So $K(\sqrt[n]{a}) = K(\sqrt[n]{a}/c)\cong K(\sqrt[n]{b})$.

    For the other direction, suppose $K(\sqrt[n]{a}) \cong K(\sqrt[n]{b})$. Then $K(\sqrt[n]{a})$ contains an $n$th root of $b$, and, by taking quotients, some $c$ an $n$th root of $a/b$. Since $v(a/b)=0$, $v(c)=0$, so $K(c)$ is unramified. As $K(c) \subseteq K(\sqrt[n]{a})$ and the latter is totally ramified, $K(c)$ it must be trivial, so $c \in K$. Then $a/b$ an $n$th power in $\OO_K$.
\end{proof}

\begin{proposition}\label{prop:k-eis-discriminant-calculation}
    Let $h(x) = \sum_{j=0}^n a_j x^j \in \OO_K[x]$ a $k$-Eisenstein polynomial, $\Delta$ its discriminant, $(k,n) = 1$. Then $v(\Delta) = k(n-1)$ and
    \[
    u = (-1)^{n(n-1)/2} n^n a_n^{n-1} w^{n-1} \mod p
    \]
    where $u$ is the unit part of $\Delta$ and $w$ the unit part of $a_0$.
\end{proposition}
\begin{proof}
    By $k$-Eisenstein, $v(a_j) \geq k - jk/n = k(n-j)/n$ for every $j$, with exact equality when $j \in \{0,n\}$. Furthermore, this condition is also necessary for equality, since $v(a_j) \in \Z$ and $k(n-j)/n \notin \Z$ for $1 \leq j < n$.

    Now, $\Delta$ is quasi-homogenous of degree $n(n-1)$ in the coefficients $a_j$, with weights $n-j$. In other words, $\Delta$ can be written as a polynomial in these coefficients, in which every term is of the form $c\prod_j a_j^{r_j}$ where $\sum_j r_j(n-j) = n(n-1)$. Then
    \[
        v\left( \prod_j a_j^{r_j} \right) = \sum_j r_j v(a_j) \geq \sum_j r_j k \dfrac{n-j}{n} = \dfrac{k}{n} \sum_j r_j(n-j) = k(n-1)
    \]
    This inequality is an equality iff $r_j = 0$ for all $j$ not 0 or $n$; since the weight of $a_n$ is 0, this is equivalent to $r_0 = n-1$.

    The actual polynomial can be computed as $a_n \Delta = {(-1)^{n(n-1)/2}} \operatorname{Res}(h', h)$,
    where $S_{f,g}$ denotes the Sylvester matrix of the polynomials $f$ and $g$, and $\operatorname{Res}(f,g) = \det S_{f,g}$ the resultant.

    Every entry of $S_{h',h}$ is either $0$, $a_j$ for $0 \leq j \leq n$, or $ja_{j-1}$ for $1 \leq j \leq n$. In particular, for every term $c\prod_j a_j^{r_j}$ of $\Delta$, the coefficient $c$ is an element of $\Z$ and thus has nonnegative valuation. And since $v(a_n) = 0$, dividing by $a_n$ has no effect on valuation. Thus
    \[v\left( c\prod_j a_j^{r_j} \right) \geq k(n-1)\]
    for all of the terms in $\Delta$, with strict inequality if $r_j \neq 0$ for some $j$ not 0 or $n$.

    The diagonal of $S_{h',h}$ consists of $n$ copies of $na_n$ followed by $n-1$ copies of $a_0$, so the corresponding term of $\operatorname{Res}(h', h)$ (in the permutation decomposition of the determinant) is $n^n a_n^n a_0^{n-1}$, which has $r_j = 0$ for all $0 < j < n$. Furthermore, no other term fulfills this condition, since all entries of $S_{h',h}$ strictly above the diagonal are either $0$ or a multiple of $a_{j}$ for some $1 \leq j < n$.

    So $\Delta$ consists of some number of terms whose valuations are $> k(n-1)$, followed by the term
    \[(-1)^{n(n-1)/2} n^n a_n^{n-1} a_0^{n-1}\]
    which has valuation exactly $k(n-1)$ since $p\nmid n$ and thus $v(n) = 0$. Therefore by the ultrametric inequality $v(\Delta) = k(n-1)$. Then multiply by $p^{-v(\Delta)} = p^{-k(n-1)}$ and recall that $v(a_0) = k$ by $k$-Eisenstein to get
    \[
    u = (-1)^{n(n-1)/2} n^n a_n^{n-1} w^{n-1} \mod p
    \]
\end{proof}
\begin{corollary}\label{cor:k-eis-shift-discriminant-properties}
    Let $f(x) \in K[x]$ irreducible with degree $n$, $\theta$ a root of $f$. If there exist $r \in \Z$, $t \in K$, $1\leq k < n$ are such that $(k,n) = 1$ and the minimal polynomial $h(x) = x^n + \sum_{i=1}^{n-1} a_i x^i$ of $p^r(\theta + t)$ is $k$-Eisenstein, then $v(\Delta_f) = k(n-1) - rn(n-1)$ (so in particular $k \equiv \frac{v(\Delta_f)}{n-1} \bmod n$), and once again
    \[
    u = (-1)^{n(n-1)/2} n^n w^{n-1} \mod p
    \]
    where $u$ is the unit part of $\Delta_f$ and $w$ the unit part of $a_0$.
\end{corollary}
\begin{proof}
    Note that $p^r(\theta + t)$ is a root of the polynomial $f(p^{-r}x - t)$, which must be irreducible since this element generates an extension of degree $n$. Linear shifts of $x$ preserve the discriminant, and scaling $x$ by some factor $c$ scales the discriminant by $c^{-n(n-1)}$, since the discriminant is homogenous in the roots with degree $n(n-1)$. So the discriminant of $f(p^{-r}x - t)$ is $p^{rn(n-1)} \Delta_f$; the result follows immediately (noting that $a_n = 1$).
\end{proof}

\begin{theorem}\label{thm:nth-power-isomorphism-criterion}
    Let $f(x), g(x) \in \OO_K[x]$ irreducible with degree $n$, $\theta$ a root of $f$, $\gamma$ a root of $g$, $K(\theta)$, $K(\gamma)$ totally ramified.

    Let $\ell$ and $s$ be inverses mod $n$ of $\frac{v(\Delta_f)}{n-1}$ and $\frac{v(\Delta_g)}{n-1}$ respectively.
    
    Then $K(\theta) \cong K(\gamma)$ (i.e. $f$ and $g$ generate isomorphic extensions of $K$) if and only if
    \[
        (-1)^{(s-\ell)\left( 1 + \frac{n(n-1)}{2}\right)} \frac{u^\ell} {w^s}
    \]
    is an $n$th power mod $p$ in $\OO_K$, where $u$ and $w$ are the unit parts of $\Delta_f$ and $\Delta_g$ respectively.
\end{theorem}
\begin{proof}
    By Proposition \ref{prop:k-eis-shift}, the minimal polynomial of $p^{r_1}(\theta + t_1)$ is $k_1$-Eisenstein for some $r_1 \in \Z$, $1 \leq k_1 < n$, $t_1 \in K$ such that $(k_1,n) = 1$. By Corollary \ref{cor:k-eis-shift-discriminant-properties}, $k_1 \equiv \frac{v(\Delta_f)}{n-1} \bmod n$, and in particular $\ell$ is also an inverse of $k_1$ mod $n$. (Since $(k_1,n)=1$, this also guarantees that such $\ell$ exists.)
    
    By Corollary \ref{cor:k-eis-nth-root-form},
    \[
    K(\theta)= K(p^{r_1}(\theta + t_1)) \cong  K\left( \sqrt[n]{p(-y_1)^\ell}\right)
    \]
    where $y_1$ is the unit part of the constant term of the minimal polynomial of $p^{r_1}(\theta + t_1)$.
    
    The same follows for $K(\gamma)$ for some $k_2$, $r_2$ and $t_2$, with $y_2$ denoting the constant term of the minimal polynomial of $p^{r_2}(\gamma + t_2)$.
    
    By Corollary \ref{cor:k-eis-shift-discriminant-properties} again,
    \[
    y_1^{n-1} \equiv (-1)^{n(n-1)/2} n^{-n} u \mod p
    \]
    so
    \[
    (-y_1)^\ell
    = (-1)^\ell(y_1^{n-1})^{-\ell} y_1^{n\ell}
    \equiv (-1)^{\ell(1+\frac{n(n-1)}{2})} n^{n\ell} u^{-\ell} \mod p
    \]
    Obviously, we get a parallel expression for $(-y_2)^s$.
    
    By Proposition \ref{prop:nth-root-extension-isomorphism-criterion}, $K(\theta) \cong K(\gamma)$ iff $(-y_2)^s/(-y_1)^\ell$ is an $n$th power mod $p$ in $\OO_K$. So we can substitute the expressions obtained above, discard $n$th powers, and simplify the quotient to see that $K(\theta) \cong K(\gamma)$ iff
     \[
        (-1)^{(s-\ell)\left( 1 + \frac{n(n-1)}{2}\right)} \frac{u^\ell} {w^s}
    \]
    is an $n$th power mod $p$.
\end{proof}

\begin{theorem}\label{thm:total-ramified-classification}
    Let $f(x) \in \OO_K[x]$ irreducible with degree $n$ and $\theta$ a root of $f$, such that $K(\theta)$ totally ramified. Let $\zeta$ a primitive $(p^m - 1)$th root of unity, where $m$ is the degree of $K/\Q_p$. Let $d = \gcd(n, p^m-1)$.
    
    If $4\nmid n$, then $K(\theta) \cong K(\sqrt[n]{\zeta^r p})$, for the unique value of $r$ mod $d$ such that
    \[
        \begin{cases}
            \zeta^{r} u^\ell & 4\nmid n\\
            -\zeta^{r} u^\ell & 4\mid n
        \end{cases}
    \]
    is an $n$th power mod $p$, where $\ell$ is an inverse mod $n$ of $\frac{v(\Delta_f)}{n-1}$, and $u$ is the unit part of $\Delta_f$.
\end{theorem}
\begin{proof}
    By Proposition 2.2.1 from \cite{Database}, all such totally ramified extensions are isomorphic to one of the form $K(\sqrt[n]{\zeta^r p})$, distinguished exactly by the value of $r$ mod $d$.

    So let $g(x) = x^n - \zeta^r p$. We can calculate directly the discriminant
    \[
        \Delta_g = (-1)^{n(n-1)/2} n^n (- \zeta^r p)^{n-1}
        = (-1)^{\frac{(n+2)(n-1)}{2}} n^n \zeta^{r(n-1)} p^{n-1}
    \]
    Take the unit part and discard $n$th powers and to obtain $(-1)^{\frac{(n+2)(n-1)}{2}} \zeta^{-r}$.

    Then Theorem \ref{thm:nth-power-isomorphism-criterion} says that $K(\theta) \cong K(\sqrt[n]{\zeta^r p})$ for $r$ (unique mod $d$) such that
    \[
        (-1)^{(1-\ell)\left( 1 + \frac{n(n-1)}{2}\right) + \frac{(n+2)(n-1)}{2}} u^\ell \zeta^{r}
    \]
    is an $n$th power mod $p$.

    When $n \equiv \{1,3\} \bmod 4$, $n$ is odd, so $-1$ is an $n$th power mod any $p$. Therefore we can ignore sign here, taking the expression to be simply $u^\ell \zeta^{r}$.

    When $n\equiv 2 \bmod 4$, $\frac{(n+2)(n-1)}{2}$ is even and $\frac{n(n-1)}{2}$ is odd, so $(-1)^{(1-\ell)\left( 1 + \frac{n(n-1)}{2}\right)} = 1$ and we again get $u^\ell \zeta^r$.

    Finally, when $n\equiv 0 \bmod 4$, $\frac{(n+2)(n-1)}{2}$ is odd, and $\ell$ is an inverse mod $n$, thus coprime to $n$, thus odd, since $n$ even. So $\ell-1$ even and $(-1)^{(1-\ell)\left( 1 + \frac{n(n-1)}{2}\right)} = -1$. Thus in this case the relevant expression changes sign to $-u^\ell \zeta^r$.
\end{proof}

\subsection{Complete classification for prime degree}

All extensions of prime degree are either unramified or totally ramified, simplifying classification. 
\footnote{Tamely ramified extensions of composite degree decompose into a unique maximal unramified subextension and a totally ramified extension above that. As such, they are likely amenable to classification via an extension of the results provided in this section; however, the problem of extracting separate generating polynomials for the ramified and unramified parts is outside the scope of this paper.}

\begin{proposition}
    For any $m \geq 1$, there is exactly one unramified extension of $\Q_p$ with degree $m$.
\end{proposition}
This is a basic fact and can be found in \cite{Neukirch} or other sources.
\begin{corollary}
    For any $m, n \geq 1$, if $K$ is an unramified degree $m$ extension of $\Q_p$, there is exactly one unramified degree $n$ extension of $K$.
\end{corollary}

\begin{theorem}
    Let $q$ prime. Let $f(x) \in \OO_K[x]$ a degree $q$ irreducible polynomial, $\theta$ a root.
    
    Then $K(\theta)/K$ is unramified if and only if $q(q-1) \mid v(\Delta_f)$.
\end{theorem}
\begin{proof}
    Suppose $K(\theta)/K$ unramified. Then it is Galois with $\Gal(K(\theta)/K) \cong C_q$, so we can pick an ordering $\theta_1, \ldots, \theta_q$ of the roots of $f$ such that $\sigma = (\theta_1, \ldots, \theta_q)$ is a generator of $\Gal(K(\theta)/K)$.
    
    Then $\sigma$ takes $(\theta_i - \theta_{i+j})$ to $(\theta_{i+1} - \theta_{i+1+j})$, where indices are mod $q$. Since an automorphism must preserve valuation, $v(\theta_i - \theta_{i+j})$ is constant in $i$, so we can denote $v(\theta_i - \theta_{i+j}) = m_j$. Furthermore,
    \[m_{j+k} = v(\theta_i - \theta_{i+j+k}) = v((\theta_i - \theta_{i+j}) + (\theta_{i+j} - \theta_{i+j+k}) \geq \min(v(\theta_i - \theta_{i+j}), v(\theta_{i+j} - \theta_{i+j+k})) = \min(m_j, m_k)\]
    so in particular $m_{sj} \geq m_j$ for any $s \in \N$.

    Since all these indices are measured mod $q$, and $q$ is prime, this implies that $m_j$ is constant in $j$. Therefore
    \begin{align*}
        \Delta_f &= \prod_{j < k} (\theta_j - \theta_k)^2 = \prod_{j=1}^{q-1} \prod_{i=1}^{q-j} (\theta_i - \theta_{i+j})^2\\
        v(\Delta_f) &= \sum_{j=1}^{q-1} \sum_{i=1}^{q-j} 2m_j = \sum_{j=1}^{q-1} 2m_j(q-j) = 2m_1 \sum_{j=1}^{q-1} j = m_1 q(q-1)
    \end{align*}
    Thus $v(\Delta_f)$ is divisible by $q(q-1)$.
    
    Conversely, if $f$ is not unramified then it is totally unramified (since $q$ prime), and then by \ref{prop:k-eis-shift} and \ref{cor:k-eis-shift-discriminant-properties} $\frac{v(\Delta_f)}{q-1} = k - rq$ where $(k,q) = 1$, so $q(q-1)\nmid v(\Delta_f)$.
\end{proof}

Combining this with the $4\nmid n$ case of Theorem \ref{thm:total-ramified-classification} yields a final result: in the case where $q$ prime and not equal to $p$, we can classify degree $q$ extensions up to isomorphism using only the discriminant of an arbitrary irreducible polynomial whose root generates the extension.

\begin{theorem}\label{thm:tame-prime-classification}
    Let $p,q$ distinct primes, $K/\Q_p$ unramified of degree $m$. Let $\zeta \in K$ a primitive $(p^m - 1)$th root of unity, and $d = \gcd(q, p^m-1)$.
    
    Suppose $f(x)\in \OO_K[x]$ an irreducible polynomial of degree $q$, $\theta \in \overline{K}$ a root of $f$.

    If $q(q-1) \mid v(\Delta_f)$, $K(\theta)/K$ is the unique unramified extension of degree $q$ over $K$.
    
    Otherwise, $K(\theta) \cong K(\sqrt[q]{\zeta^r p})$ for the (unique mod $d$) value of $r$ such that $\zeta^r u^\ell$
    is an $q$th power mod $p$, where $\ell$ is the inverse mod $q$ of $\frac{v(\Delta_f)}{q-1}$ and $u$ is the unit part of $\Delta_f$.
\end{theorem}

\section{Toward cubic extensions of $\Q_3$}

From this section onward, we consider the degree $3$ extensions of $\Q_3$. Our methods using generalized Eisenstein polynomials does not apply in the presence of wild ramification, so in order to classify these extensions, we look for algebraic invariants.

\subsection{Isomorphism classes}
We first adapt the classification result proved in \cite{Amano} for degree 3:
\begin{proposition}\label{prop:wildly ramified classification table}
There are 9 (wildly, totally) ramified cubic extensions of $\Q_3$ up to isomorphism:
\begin{center}
\begin{tabular}{|c|c|c|c|}
\hline
Polynomial       & Ramification Exponent & Galois Group & Inertia Group \\ \hline
$x^3 + 3x + 3$   & $3$                     & $S_3$        & $S_3$         \\ \hline
$x^3 + 6x + 3$   & $3$                      & $S_3$        & $S_3$         \\ \hline
$x^3 + 3x^2 + 3$ & $4$                     & $S_3$        & $C_3$         \\ \hline
$x^3-3x^2+3$     & $4$                     & $C_3$        & $C_3$         \\ \hline
$x^3-3x^2+12$    & $4$                     & $C_3$        & $C_3$         \\ \hline
$x^3-3x^2+21$    & $4$                     & $C_3$        & $C_3$         \\ \hline
$x^3+3$          & $5$                     & $S_3$        & $S_3$         \\ \hline
$x^3+12$         & $5$                     & $S_3$        & $S_3$         \\ \hline
$x^3+21$         & $5$                     & $S_3$        & $S_3$             \\ \hline
\end{tabular}
\end{center}
\end{proposition} 
\medskip
Note that there are 3 Galois extensions, which can be distinguished from the rest by their discriminant:

\begin{lemma} \label{Splitting Field of a Cubic}
Let \(f(x)\) be a separable cubic in \(F[x]\) with a root \(\theta\) and discriminant \(\Delta\). Then the splitting field of \(f\) over \(F\) is \(F(\theta, \sqrt{\Delta})\), and $\Delta$ is a square in $F$ if and only if $F(\theta)$ is Galois. 
\end{lemma}
\begin{proof}
    This is a standard result in field theory. For proof, see for example \cite{Conrad}.
\end{proof}

\begin{corollary}
When $\Delta_f$ is not a square, we have the following arrangement of extensions of $\Q_3$, and we call $\Q_3(\sqrt{\Delta})$ the degree 2 subextension (under the Galois closure):
    
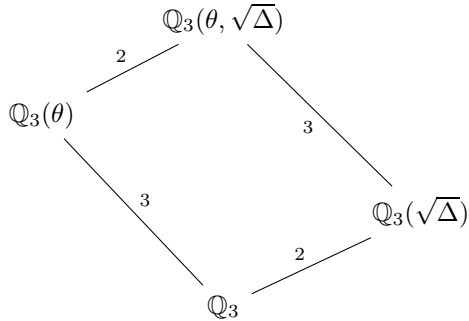
\[\begin{tikzcd}
	& {\Q_3(\theta, \sqrt{\Delta})} \\
	{\Q_3(\theta)} & {} \\
	&& {\Q_3(\sqrt{\Delta})} \\
	& \Q_3
	\arrow["3", no head, from=2-1, to=4-2]
	\arrow["2 ", no head, from=4-2, to=3-3]
	\arrow["3", no head, from=3-3, to=1-2]
	\arrow["2", no head, from=2-1, to=1-2]
\end{tikzcd}\]
\end{corollary}

Note that the quadratic subextension $\Q_3(\Delta_f)$ is an invariant of the field extension $\Q_3(\theta_f)$: if $\Q_3(\Delta_g) \ncong \Q_3(\Delta_f)$ then $\Q_3(\theta_g) \ncong \Q_3(\theta_f)$. The classification of quadratic extensions of $\Q_p$ is also well-known and easy to determine just in terms of the discriminant, so we can use quadratic subextension as a tool to determine isomorphism classes.

\begin{proposition} \label{cubicTable}
    The 9 ramified extensions of $\Q_3$ up to isomorphism can be grouped in the following way:
\begin{center}
\begin{tabular}{ |c|c|c|c| } 
\hline
Galois & Quadratic subextension & Generating polynomial \\
\hline
\multirow{3}{1.3em}{Yes} & \multirow{3}{1.2em}{N/A} & $x^3 - 3x^2 +3$ \\ 
& & $x^3 - 3x^2 + 12$  \\ 
& & $x^3 - 3x^2 + 21$ \\ 
\hline
No & $\Q_3(\sqrt{-1})$ & $x^3 + 3x^2 +3$ \\ 
\hline
No & $\Q_3(\sqrt{3})$ & $x^3 + 6x +3$ \\ 
\hline
\multirow{5}{1.2em}{No}& \multirow{5}{3em}{$\Q_3(\sqrt{-3}$)} & $x^3 + 3x +3$ \\ 
& & $x^3 + 3$ \\ 
& & $x^3 + 12$ \\ 
& & $x^3 + 21$ \\ 
\hline
\end{tabular}
\end{center}
\end{proposition}

We can immediately distinguish two of the classes:

\begin{corollary} \label{Wild Classify Thm 1}
 Let $f(x)$ be a degree 3 irreducible polynomial in $\Q_3[x]$. Denote by $\theta$ a root of $f$ and by $\Delta$ its discriminant.
 
 $\Delta$ is not a square in $\Q_3$ iff one the following is true:
\begin{enumerate}
    \item $v_3(\Delta)$ odd and $\Q_3(\theta) \cong \Q_3(x) / (x^3 + 6x + 3)$
    \item $v_3(\Delta)$ even and $\Q_3(\theta) \cong \Q_3(x) / (x^3 + 3x^2 + 3)$
\end{enumerate}
\label{corollary: classification of x^3 + 3x^2 + 3}
\end{corollary}

\subsection{Preliminaries}

To distinguish between the other isomorphism classes, we recall some basic results from class field theory. 

\begin{theorem} (Local Artin reciprocity)
    Let $F$ be a local field. There is a unique continuous homomorphism
    \begin{equation*}
        \Theta_F: F^{\times} \rightarrow \Gal(F^{ab}/F)
        \label{eq: homomorphism for artin reciprocity}
    \end{equation*}
    such that for each finite Galois extension $L/F$ in $F^{\mathrm{ab}}$, the homomorphism
    \begin{equation*}
        \Theta_{L/F}: F^{\times} \rightarrow \Gal(L/F)
    \end{equation*}
    obtained by composing $\Theta_F$ with quotient map $\Gal(F^{ab}/F) \rightarrow \Gal(L/F)$ satisfies: 
    \begin{enumerate}
        \item if $F$ is nonarchimedean and $L/F$ is unramified, then $\Theta_{L/F} (\pi) = \operatorname{Frob}_{L/F}$ for every uniformizer $\pi$ of $\OO_F$
        \item $\Theta_{L/F}$ is surjective with kernel $N_{L/F}(L^{\times})$, inducing $F^{\times}/N_{L/F}(L^{\times})\simeq \Gal(L/F)$
    \end{enumerate}
    \label{Theorem: Local Artin Reciprocity}
\end{theorem}

\begin{corollary} \label{norm group cube lemma}
If the degree of $L/F$ is $n$,
\[
\frac{F^{\times} / (F^{\times})^n}{N_{L/F}(L^{\times}) / (F^{\times})^n} \cong \Gal(L/F)
\]
\end{corollary}

\begin{proof}
For each $a \in F$, $N_{L/F}(a) = a^n$, so 
$(F^\times)^n$ is a subgroup of $N_{L/F}(L^\times)$. Then the result follows from the isomorphism theorems and part 2 of Theorem \ref{Theorem: Local Artin Reciprocity}.
\end{proof}

\begin{proposition}\label{prop:norm group uniqueness}
    Let $F$ a local field. Finite abelian extensions $L/F$ within $F^{\mathrm{ab}}/F$ are in inclusion-reversing bijection with the norm groups of $F$ via
    \[
    L/F \leftrightarrow N_{L/F}(L^\times)\leq F^\times
    \]
\end{proposition}
\begin{proof}
    This is a standard consequence of the Norm Limitation Theorem, which can be found for example in \cite{Neukirch}.
\end{proof}

We now recall the following properties of $p$-adic fields. For proof, see for example \cite{Neukirch} Chapter 2 Section 5. 

\begin{proposition} \label{DecompLemma}
    Let $F$ be a $p$-adic field with residue field $\F_q$, valuation ring $\OO_F$ and maximal ideal $\p$. Let $p\OO_F = \p^e$. Then for $n > \frac{e}{p-1}$,
    \[
    \exp:\p^n \to U^{(n)}
    \]
    is an isomorphism, where $U^{(n)} = 1 + \p^n$ with addition. Also,
    \[
    \p^n = \pi^n \OO_F \cong \OO_F \cong \Z_p^m,
    \]
    where $m = [F:\Q_p]$. 
\end{proposition}

\begin{proposition} \label{Decomposition of p-adic field}
    Under the same hypotheses,
    \[
    F^\times \cong \Z \oplus \Z / (q-1)\Z \oplus U^{(1)}.
    \]
\end{proposition}
From this following two useful decompositions: 
\begin{corollary}\label{decompQp}
    For $p>2$, we can can decompose the group $\Q_p^\times$ as 
    \[
    \Q_p^\times \cong \Z \oplus \Z/(p-1)\Z \oplus \Z_p,
    \]
    and when $p=2$,
    \[
    \Q_2^\times \cong \Z \oplus \Z/2\Z \oplus \Z_2.
    \]
\end{corollary}

\begin{corollary} \label{decomp Q_p(sqrt)} 
    For $p>3$, we can decompose $\Q_p(\sqrt{cp})^\times$ as
    \[
    \Q_p(\sqrt{cp})^\times \cong \Z \oplus \Z/(p-1)\Z \oplus \Z_p^2.
    \]
    and when $p = 3$,
    \[
    \Q_3(\sqrt{3c})^\times \cong \Z \oplus \Z/2\Z \oplus \Z / 3\Z \oplus \Z_3^2
    \]
    
\end{corollary}

\begin{proof}
    For $p > 3$, by Prop. $\ref{DecompLemma}$, $U^{(1)} \cong \p^2 \cong \Z_p^2$.
    
    For $p = 3$, $U^{(2)} \cong \p^2 \cong \Z_p^2$. Since $U^{(1)} / U^{(2)} \cong \OO_K / \p = \Z/3\Z$, $U^{(1)} \cong \Z/3\Z \oplus \Z_p^2$.
\end{proof}

\medskip
Using these, we obtain the structure of cube class groups in two useful cases:
\begin{lemma}\label{cubeclass}
    $\Q_3^\times / (\Q_3^{\times})^3 \cong C_3 \oplus C_3$.
\end{lemma}
\begin{proof}
    From \ref{decompQp},
    \[
    \Q_3^\times \cong \Z \oplus \Z_3 \oplus \Z/2\Z,
    \]
    so
    \[
    (\Q_3^\times)^3 \cong 3\Z \oplus 3\Z_3 \oplus 3\Z/2\Z = 3\Z \oplus 3\Z_3 \oplus \Z/2\Z.
    \]
    \[
    \Q_3^\times / (\Q_3^{\times})^3 \cong \Z / 3\Z \oplus \Z_3 / 3\Z_3 = C_3 \oplus C_3
    \]
\end{proof} 

\begin{lemma} \label{sqrt-3 decomp}
      $\Q_3 (\sqrt{3c})^\times / (\Q_3^{\times}(\sqrt{3c}))^3 \cong C_3 \oplus  C_3 \oplus C_3 \oplus C_3$.
\end{lemma}
\begin{proof}
    By the decomposition in \ref{decomp Q_p(sqrt)},
    \[
    \Q_3 (\sqrt{3c})^\times \cong \Z \oplus \Z/2\Z \oplus \Z/3\Z \oplus \Z_3 \oplus \Z_3,
    \]
    so
    \[
    (\Q_3 (\sqrt{3c})^\times)^3 \cong 3\Z \oplus 3\Z/2\Z \oplus 3\Z/3\Z \oplus 3\Z_3 \oplus 3\Z_3 = 3\Z \oplus \Z/2\Z \oplus 1 \oplus 3\Z_3 \oplus 3\Z_3.
    \] 
    \[
    \Q_3 (\sqrt{3c})^\times / (\Q_3^{\times}(\sqrt{3c}))^3 \cong C_3 \oplus C_3 \oplus C_3 \oplus C_3.
    \]
\end{proof}

Since depressed forms of polynomials are easy to compute and convenient, we will use them heavily. We introduce a basic lemma on the coefficients of the polynomial in the case where it is depressed:
\begin{lemma} \label{lem:wild-newton-polygon}
    Let $f(x) = x^3 + \alpha x + \beta$ an irreducible polynomial over $\Q_3$ that generates a totally ramified extension. Then $v_3(\alpha) - \frac{2}{3}v_3(\beta) > 0$. Here $v_3(\alpha)$ may be $\infty$, if $\alpha = 0$.
\end{lemma}
\begin{proof}
    By irreducibility, the Newton polygon must be a single line segment; otherwise it would have a segment of length 1, and thus a root. Hence $v_3(\alpha) - \frac{2}{3}v_3(\beta) \geq 0$.

    Now suppose that $v_3(\alpha) = \frac{2}{3}v_3(\beta)$.Then $v_3(\alpha) = 2m$, $v_3(\beta) = 3m$ for some $m \in \Z$. Write $\alpha = 3^{2m}w$, $\beta = 3^{3m}u$ for $w,u$ units.
    
    If $\theta$ is a solution to $f(x) = x^3 + \alpha x + \beta$, then $3^{-m} \theta$ is a solution to $3^{3m} x^3 + 3^m \alpha x + \beta$ and thus to $g(x) = x^3 + 3^{-2m} \alpha x + 3^{-3m}\beta = x^3 + wx + u$, and they generate the same extension. So it suffices to consider the case where $f(x) = x^3 + wx + u$.
    
    In particular, this means $f(x) \in \Z_3[x]$ and $f'(x) = 3x^2 + w \equiv w \not \equiv 0$ mod 3, so $f(x)$ is separable over $\F_3$. Then since $f(x)$ is irreducible over $\Q_3$, it is irreducible over $\F_3$ by Hensel's Lemma, and thus $f(x)$ generates an unramified extension over $\Q_3$, a contradiction.
\end{proof}

\section{Galois cases}
Throughout this section, $\Q_3(\theta)/\Q_3$ is a ramified Galois extension, and $f(x) \in \Q_3[x]$ is an irreducible cubic polynomial with root $\theta$, unless otherwise stated. (Recall that by Lemma \ref{Splitting Field of a Cubic} this is equivalent to $\Delta_f$ being a square.)
 
We also denote the norm group $N_{\Q_3(\theta) / \Q_3} (\Q_3(\theta)^\times) \leq \Q_3^\times$ as simply $N(\Q_3(\theta)^\times)$.

\subsection{The cubic norm group}
    
\begin{definition}
    We refer to
    \[
    N(\Q_3(\theta)^\times) / (\Q_3^\times)^3 < \Q_3^\times/(\Q_3^\times)^3
    \]
    as the \textbf{cubic norm group} of the extension $\Q_3(\theta)/\Q_3$. When the extension is generated by an irreducible cubic $f$, we may also call it the cubic norm group associated with $f$.
\end{definition}

\begin{lemma}\label{lem:cubic norm group size lemma}
\[
N(\Q_3(\theta)^\times) / (\Q_3^\times)^3 \cong C_3.
\]
\end{lemma}
\begin{proof} By Corollary \ref{norm group cube lemma}, 
\[
\frac{\Q_3^\times / (\Q_3^\times)^3}{N(\Q_3(\theta)^\times) / (\Q_3^\times)^3 } \cong \Gal (\Q_3(\theta) / \Q_3)
\]
The extension of $\Q_3(\theta) / \Q_3$ is Galois and is of degree 3, so $\Gal (\Q_3(\theta) / \Q_3) \cong C_3$. Therefore $N(\Q_3(\theta)^\times) / (\Q_3^\times)^3$ is a subgroup of index 3 in $\Q_3^\times / (\Q_3^\times)^3$. By Corollary \ref{decompQp}, $\Q_3^\times / (\Q_3^{\times})^3 \cong C_3 \oplus C_3$, and the result follows.
\end{proof}

\begin{proposition}
    If $\Q_3(\theta_1)/\Q_3$ and $\Q_3(\theta_2)/\Q_3$ are nonisomorphic, their cubic norm groups intersect trivially in $\Q_3^\times/(\Q_3^\times)^3$.
\end{proposition}
\begin{proof}
    Since the cubic norm group of each extension is $\cong C_3$, which has no nontrivial subgroups, it suffices to show that they are distinct. But this follows from Proposition \ref{prop:norm group uniqueness}; if the norm groups are distinct, their quotients by $(\Q_3^\times)^3$ must be also.
\end{proof}

\begin{lemma} \label{cubes in Q3}
     $a \in \Z_3^\times$ is a cube iff $a \equiv \pm 1 \bmod 9$. 
     \label{cubes in Q3}
\end{lemma}
\begin{proof}
    One direction is obvious. For the other, take $a \equiv \pm 1 \bmod 9$ and let $f(x) = x^3 - a$. There exists $r$ with $f(r) = r^3 - a \equiv 0 \bmod 9$; since $r \not\equiv 0 \bmod 3$, $f'(r) = 3r^2 \not\equiv 0 \bmod 9$. Then by Hensel's lemma $a$ is a cube.
\end{proof}

\begin{theorem} \label{thm:galois cubic norm group computation}
    Let $\theta$ a root of $x^3 - 3x^2 + 3\tau$, where
    $\tau \in \{1, 4, 7\}$. Then the three 
    elements of 
    $N(\Q_3(\theta)^\times)/(\Q_3^\times)^3$ 
    have representatives $\{1, 3\tau, 9 \tau^2\}$.
\end{theorem}
\begin{proof}
From \ref{lem:cubic norm group size lemma}, $N(\Q_3(\theta)^\times)/(\Q_3^\times)^3 \cong C_3$, so it suffices to show that $3\tau$ represents a nontrivial element; then it is a generator and $9 \tau^2$ distinct and nontrivial.

We compute $N(\theta) = -3 \tau$. Since $v_3(\tau) = 0$, $-3\tau$ has valuation $1$ and thus cannot be a cube, so it is nontrivial in $\Q_3^\times/(\Q_3^\times)^3$. Since $-1$ a cube, the same is true of $3\tau$.
\end{proof}

\subsection{Determining the isomorphism class}
Given an arbitrary depressed cubic generating a Galois extension, we can now determine that extension up to isomorphism by finding one nontrivial element of its cubic norm group and comparing it with the ones we just computed.

\begin{lemma} \label{r-2m = 1 lemma}
    Let $f(x) = x^3 + \alpha x + \beta$.
    Then if $v_3(\beta) \equiv 0$ mod 3, $v_3(\alpha) - \frac{2}{3}v_3(\beta) = 1$. 
\end{lemma}
\begin{proof}
    Suppose $v_3(\beta) \equiv 0$ mod 3. Then $v_3(\alpha) - \frac{2}{3}v_3(\beta)$ is an integer. By Lemma \ref{lem:wild-newton-polygon}, $v_3(\alpha) - \frac{2}{3}v_3(\beta) > 0$.
    Then, $\alpha = 3^r w$ and $\beta = 3^{3m}u$, where $u,w$ are 3-adic units and $r-2m \geq 2$. ($r$ may be $\infty$.)
    
    Suppose $r - 2m \geq 2$. The discriminant of $f(x)$ is
    \[
    \Delta = -4 \alpha^3 - 27 \beta^2 = -4 \cdot 3^{3r} \cdot w^3 - 3^{6m+3}u^2 = 3^{6m+3} (-4w^3 \cdot 3^{3(r-2m-1)} - u^2).
    \]
    where $v_3(-4w^3 \cdot 3^{3(r-2m-1)}) > 0$ (since $r-2m \geq 2$) and $v_3(u^2) = 0$. So $v_3(-4w^3 \cdot 3^{3(r-2m-1)} - u^2) = 0$ and thus $v_3(\Delta) = 6m+3$, an odd number, which is a contradiction, since $\Delta$ must be square in $\Q_3$.
    
    If $\alpha = 0$, $\Delta = -27\beta^2 = -3^{6m + 3}u^2$ and once again $v_3(\Delta) = 6m+3$, so the same argument applies.
\end{proof}

\begin{proposition} \label{Galois 3 cases enumeration}
    Let $f(x) = x^3 + \alpha x + \beta$, where $w$ and $u$ are the unit parts of $\alpha$ and $\beta$ respectively.
    
\begin{enumerate}
    \item If $v_3(\beta) \not \equiv 0$ mod 3, then $\beta$ represents a nontrivial element of the cubic norm group of $\Q_3(\theta)/\Q_3$.

    The same element is represented equivalently by $3^\ell u$ for the $\ell \in \{1,2\}$ such that $v_3(\beta) = 3m + \ell$.
    
    \item If $v_3(\beta) \equiv 0$ mod 3 and $u \not \equiv \pm 1$ mod 9, then $\beta$ (equivalently, $u$) represents a nontrivial element of the cubic norm group of $\Q_3(\theta)/\Q_3$.
    
    \item If $v_3(\beta) \equiv 0$ mod 3 and $u \equiv 1$ mod 9, then $1 - u + 3w$ represents a nontrivial element of the cubic norm group of $\Q_3(\theta)/\Q_3$, and has valuation 1.
    
    \item If $v_3(\beta) \equiv 0$ mod 3 and $u \equiv -1$ mod 9, then $1 + u + 3w$ represents a nontrivial element of the cubic norm group of $\Q_3(\theta)/\Q_3$, and has valuation 1.
\end{enumerate}
(In particular, in cases 3 and 4, $v_3(\alpha) = \frac{2}{3}v(\beta) + 1 < \infty$, so the unit part of $\alpha$ exists.)
\end{proposition}

\begin{proof}
\begin{enumerate}
    \item By direct calculation, $N(\theta) = -\beta$. Since $3\nmid v_3(\beta) = v_3(-\beta)$, $-\beta$ is not a cube. Factoring out the cube $3^{3m}$ yields $-3^\ell u$, and since $-1$ a cube we can ignore sign.

    \item Again, $N(\theta) = -\beta$, which differs by the cube $-3^{3m}$ from $u$. Since $u \not \equiv \pm 1$ mod 9, $u$ is not a cube in $\Z_3^\times$, so $u,\beta$ represent a nontrivial element.
    
    \item We can compute $N(\theta + \lambda) = -\beta + \lambda^3 + \lambda\alpha$. So (multiplying by $3^{-3m}$, which is cube),
    \[
    3^{-3m}N(\theta + 3^m)
    = 3^{-3m}(-\beta + 3^{3m} + 3^m\alpha)
    = -u + 1 + 3^{-2m}\alpha
    \]
    By Lemma \ref{r-2m = 1 lemma}, $v_3(\alpha) - \tfrac{2}{3}v_3(\beta) = 1$, so this is equal to $-u + 1 + 3w$.
    
    Since $u \equiv 1 $ mod 9 and $w$ is a unit, $3w - u + 1 \equiv 3w$ mod 9. Since $w$ a unit, $3w-u+1$ has valuation 1 and is not a cube.
    
    \item Similarly to the previous case,
    \[
    -3^{-3m}N(\theta - 3^m)
    = -3^{-3m}(-\beta - 3^{3m} - 3^m\alpha)
    = u + 1 + 3w
    \]
    again using Lemma \ref{r-2m = 1 lemma}. Since $u \equiv -1$ mod 9 and $w$ is a unit, $u + 1 + 3w\equiv 3w$ mod 9, so it has valuation 1 and is not a cube.
\end{enumerate}
\end{proof}

Case 2 of the above proposition is actually impossible:

\begin{corollary}
    Let $f(x) = x^3 + \alpha x + \beta$ be an irreducible polynomial over $\Q_3$, $\theta\in\bar\Q_3$ a root of $f$, $u$ the unit part of $\beta$. If $v_3(\beta) \equiv 0$ mod 3 and $u \not \equiv \pm 1$ mod 9, then $f(x)$ does not generate a Galois, wildly ramified extension.
\end{corollary}
\begin{proof}
    By Lemma \ref{r-2m = 1 lemma}, $v_3(\alpha) - \frac{2}{3}v_3(\beta) = 1$, so we can write $\alpha = 3^{2m+1}w$ and $\beta = 3^{3m}u$.

    Suppose that $\theta$ generates a Galois, wildly ramified extension $\Q_3(\theta)/\Q_3$. Then it is isomorphic to one of the three Galois cases, so by Theorem \ref{thm:galois cubic norm group computation},
    \[
    N(\Q_3(\theta)) / (\Q_3^\times)^3 = \{1, 3\tau, 9 \tau^2 \}(\Q_3^\times)^3
    \]
    for some $\tau \in \{1, 4, 7\}$. 
    
    By Proposition \ref{Galois 3 cases enumeration}, $u$ represents a nontrivial element of the cubic norm group of $\Q_3(\theta)/\Q_3$. So it must differ by a cube from either $3\tau$ or $9\tau^2$. But this is a contradiction, since $v_3(u) = 0$, $v_3(3\tau) = 1$, $v_3(9\tau^2) = 2$, which are all different mod 3.
\end{proof}

We therefore present all possible cases such that $f(x)$ generates a Galois, wildly ramified extension:
\begin{theorem} \label{main theorem Galois}
    Let $f(x) = x^3 + \alpha x + \beta$ be a depressed irreducible cubic over $\Q_3$ with root $\theta$. Let $u$ the unit part of $\beta$, $w$ the unit part of $\alpha$ (if $\alpha\neq 0$).
    
    Suppose that $\theta$ generates a ramified Galois extension $\Q_3(\theta)/\Q_3$, or equivalently  $\Delta_f$ square in $\Q_3$. Then $3\mid v_3(\beta) \implies \alpha\neq 0$, and if we define
    \begin{align*}
        t = \begin{cases}
            w + \frac{1-u}{3} & v_3(\beta) \equiv 0 \bmod 3,\quad u \equiv 1 \bmod 9\\
            w + \frac{1+u}{3} & v_3(\beta) \equiv 0 \bmod 3,\quad u \equiv -1 \bmod 9\\
            u & v_3(\beta) \equiv 1 \bmod 3\\
            u^2 & v_3(\beta) \equiv 2 \bmod 3
        \end{cases}
    \end{align*}
    then $\Q_3(\theta)/\Q_3 \cong \Q_3[x]/(x^3-3x^2+ \beta)$ for the unique $\beta \in \{1, 4, 7\}$ such that $t \equiv \pm \tau \bmod 9$. 
\end{theorem}

\begin{proof}
    Suppose $v_3(\beta) \equiv 0 \bmod 3$. Then by Proposition \ref{Galois 3 cases enumeration}, $u\equiv \pm 1 \bmod 9$, $\alpha \neq 0$, and $3t = 1 \mp u + 3w$ has valuation 1 and represents a nontrivial element of the cubic norm group of $\Q_3(\theta)/\Q_3$.

    If $v_3(\beta) \equiv 1 \bmod 3$, then by Proposition \ref{Galois 3 cases enumeration} case 1, the same is true of $3t = 3u$.

    Finally, if $v_3(\beta) \equiv 2 \bmod 3$, then by Proposition \ref{Galois 3 cases enumeration} case 1, $3^2 u$ is a nontrivial element of the cubic norm group, so its square $3^4 u^2$ must be as well, and thus so is $3t = 3u^2$.

    So in each case $3t$ represents a nontrivial element of the cubic norm group and has valuation 1. Since $\Q_3(\theta)/\Q_3$ is Galois, it is in the isomorphism class of the extension generated by $x^3 - 3x^2 + 3\tau$ for some $\beta \in \{1, 4, 7\}$. By Theorem \ref{thm:galois cubic norm group computation}, this means that $t$ differs by a cube from one of $\{1, 3\tau, 9\tau^2 \}$. Since $v_3(\tau) = 1$ and $v_3(3t) = 1$, it must be $3\tau$, so $3t/3\tau = t/\tau$ is a cube in $\Q_3^\times$. By Lemma \ref{cubes in Q3}, this is equivalent to $t/\tau \equiv \pm 1 \bmod 9$, or $t \equiv \pm  \tau \bmod 9$.
\end{proof}

\section{Remaining non-Galois cases}
Throughout this section, $\Q_3(\theta)/\Q_3$ is a ramified non-Galois extension whose Galois closure has quadratic subfield $\cong \Q_3(\sqrt{-3})$, and $f(x) \in \Q_3[x]$ is an irreducible cubic polynomial with root $\theta$, unless otherwise stated.

This comprises four isomorphism classes:
\begin{center}
    \begin{tabular}{|c|c|c|c|}
         \hline
        Galois & Quadratic subextension & Generating polynomial \\
        \hline
        \multirow{4}{1.2em}{No}& \multirow{4}{3em}{$\Q_3(\sqrt{-3}$)} & $x^3 + 3x +3$ \\ 
        & & $x^3 + 3$ \\ 
        & & $x^3 + 12$ \\ 
        & & $x^3 + 21$ \\ 
        \hline
    \end{tabular}
    \label{tab:non galois classes}
\end{center}

Instead of examining $\Q_3(\theta)/\Q_3$ directly, we consider the related extension $\Q_3(\theta, \sqrt{\Delta})/\Q_3(\sqrt{\Delta})$, which is Galois and abelian. This allows us to once again use local Artin reciprocity to differentiate isomorphism classes via quotients of the norm group.

\begin{lemma}\label{lem:non-galois discriminant criterion}
    An irreducible cubic $f(x) \in \Q_3[x]$ generates an extension isomorphic to that generated by $x^3 + 3x + 3$, $x^3 + 3$, $x^3 + 12$, or $x^3 + 21$ iff $\frac{\Delta_f}{-3}$ is a square in $\Q_3^\times$.
\end{lemma}
\begin{proof}
    By assumption, $\Q_3(\sqrt{\Delta})\cong\Q_3(\sqrt{-3})$, so (being Galois) they are equal as subfields of $\Q_3(\theta,\sqrt{\Delta})$.
\end{proof}

\subsection{The reduced cubic norm group}
From now on, we denote $N_{\Q_3(\theta, \sqrt{-3}) / \Q_3(\sqrt{-3})} (\Q_3(\theta, \sqrt{-3})^\times)$ as simply $N(\Q_3(\theta, \sqrt{-3})^\times)$. 

\begin{lemma}
    $N(\Q_3(\theta, \sqrt{-3})^\times)$ contains $\Q_3^\times$.
\end{lemma}
\begin{proof}
    This is a consequence of standard facts on the multiplicative structure of $p$-adic fields, which can be found for instance in Chapter 5 of \cite{Neukirch}. Note in particular that in the decomposition
    \[
    \Q_3(\sqrt{-3})^\times \cong \Z \oplus \Z/2\Z \oplus \Z / 3\Z \oplus \Z_3^2
    \]
    from Corollary \ref{decomp Q_p(sqrt)}, the factor of $\Z/3\Z$ corresponds to $U^{(1)}/U^{(2)}$, represented by the third roots of unity. Since the norm map from $\Q_3(\theta,\sqrt{-3})$ multiplies together Galois conjugates, this is exactly the part of the decomposition it doesn't cover (its image being an index 3 subgroup of $\Q_3(\sqrt{-3})^\times$) and the part that is not present in $\Q_3$.
\end{proof}

\begin{definition}
    We refer to the subgroup
    \[
    \frac{N(\Q_3(\theta, \sqrt{-3})^\times )}{\Q_3^\times (\Q_3(\sqrt{-3})^\times)^{3}} < \frac{\Q_3(\sqrt{-3})^\times}{\Q_3^\times (\Q_3(\sqrt{-3})^\times)^{3}}
    \]
    where $\Q_3^\times(\Q_3(\sqrt{-3})^\times)^{3}$ denotes the subgroup of $\Q_3(\sqrt{-3})^\times$ generated by $(\Q_3(\sqrt{-3})^\times)^3$ and $\Q_3^\times$,
    as the \textbf{reduced cubic norm group} of $\Q_3(\theta, \sqrt{-3})/\Q_3(\sqrt{-3})$.
    
    If $\theta$ is the root of an irreducible cubic polynomial $f(x) \in \Q_3[x]$, we may also call this group the reduced cubic norm group associated with $f(x)$.
\end{definition}

\begin{remark}
In other words, a nontrivial element of the reduced cubic norm group is represented by $N(x) \in \Q_3(\sqrt{-3})$ for some $x \in \Q_3(\theta, \sqrt{-3})$ such that $N(x)$ does not differ by an element of $\Q_3^\times$ from any cube in $\Q_3(\sqrt{-3})^\times$.
\end{remark}

We will show that our remaining four isomorphism classes exhibit distinct reduced cubic norm groups, and use this to identify which one the extension generated by a given polynomial belongs to.
\medskip

\begin{lemma} \label{reduced cubic norm group isom to C3}
    If $\Q_3(\theta)/\Q_3$ is cubic and non-Galois with quadratic subextension $\Q_3(\sqrt{-3})$,
    \[
    \frac{N(\Q_3(\theta, \sqrt{-3})^\times )}{\Q_3^\times(\Q_3(\sqrt{-3})^\times)^{3}}  \cong C_3
    \]
\end{lemma}
\begin{proof}
By Corollary \ref{sqrt-3 decomp}, 
\[
    \Q_3 (\sqrt{-3})^\times / (\Q_3^{\times}(\sqrt{-3}))^3 \cong C_3 \oplus C_3 \oplus  C_3 \oplus C_3.
\]
By Corollary \ref{norm group cube lemma}, 
\[
    \frac{\Q_3(\sqrt{-3})^\times / (\Q_3(\sqrt{-3})^\times)^3}{N(\Q_3(\theta, \sqrt{-3})^\times) / (\Q_3(\sqrt{-3})^\times)^3 } \cong \Gal (\Q_3(\theta, \sqrt{-3}) / \Q_3 (\sqrt{-3})) \cong C_3
\]
so $N(\Q_3(\theta, \sqrt{-3})^\times) / (\Q_3(\sqrt{-3})^\times)^3$ must be a subgroup of index 3,
and therefore
\[
    N(\Q_3 (\theta, \sqrt{-3})^\times) / (\Q_3(\sqrt{-3})^\times)^3 \cong C_3 \oplus C_3 \oplus C_3.
\]
Now, if $a + b\sqrt{-3} \in \Q_3(\sqrt{-3})^\times$, 
\[
(a + b\sqrt{-3})^3 = a^3 - 9ab^2 + (3a^2b - 3b^3)\sqrt{-3}
\]
This is in $\Q_3^\times$ iff $b=0$ or $a^2 = b^2$. If the former, it is $a^3$. If the latter, it is $a^3 - 9a^3 = 8a^3$. In both cases it is in $(\Q_3^\times)^3$. In other words, $(\Q_3(\sqrt{-3})^\times)^3 \cap \Q_3^\times = (\Q_3^\times)^{3}$. From this, it's not hard to see that 
\[
\Q_3^\times(\Q_3(\sqrt{-3})^\times)^{3}/(\Q_3(\sqrt{-3})^\times)^{3} \cong \Q_3^\times / (\Q_3^{\times})^3 
\]
and by Corollary \ref{decompQp},
\[
\Q_3^\times / (\Q_3^{\times})^3 \cong C_3 \oplus C_3,
\]
Then by the third isomorphism theorem,
\[
    \frac{N(\Q_3(\theta, \sqrt{-3})^\times )}{\Q_3^\times (\Q_3(\sqrt{-3})^\times)^{3}}
    \cong \frac{N(\Q_3(\theta, \sqrt{-3})^\times )/(\Q_3(\sqrt{-3})^\times)^{3}}{\Q_3^\times (\Q_3(\sqrt{-3})^\times)^{3}/(\Q_3(\sqrt{-3})^\times)^{3}}
\]
has order $|C_3 \oplus C_3 \oplus C_3|/|C_3 \oplus C_3| = 3$ and is thus isomorphic to $C_3$.
\end{proof}

\begin{proposition}
    If $\Q_3(\theta_1)/\Q_3$ and $\Q_3(\theta_2)/\Q_3$ are nonisomorphic, the reduced cubic norm groups of $\Q_3(\theta_1, \sqrt{-3})/\Q_3(\sqrt{-3})$ and $\Q_3(\theta_2, \sqrt{-3})/\Q_3(\sqrt{-3})$ intersect trivially in $\frac{\Q_3(\sqrt{-3})^\times}{\Q_3^\times (\Q_3(\sqrt{-3})^\times)^3}$.
\end{proposition}
\begin{proof}
    Since the reduced cubic norm groups are $\cong C_3$, it suffices to show that they are distinct. If $\Q_3(\theta_1)/\Q_3$ and $\Q_3(\theta_2)/\Q_3$ are nonisomorphic, they have different Galois closures, so $\Q_3(\theta_1, \sqrt{-3}) \neq \Q_3(\theta_2, \sqrt{-3})$. Since both are abelian extensions of $\Q_3(\sqrt{-3})$, by Proposition \ref{prop:norm group uniqueness} their norm groups must be distinct, and thus the quotients by $\Q_3^\times (\Q_3(\sqrt{-3})^\times)^3$ must also be distinct.
\end{proof}

To compute these groups effectively, we use a criterion describing cubes in $\Q_3(\sqrt{-3})$:
\begin{lemma} \label{cubes in Q_3(sqrt-3)}
    An element $c \in \Q_3(\sqrt{-3})^\times$ is a cube if and only if $c = ({\sqrt{-3}})^{3m} u$ for $u \in \Z_3[\sqrt{-3}]^{\times}$ such that $u \equiv \pm 1\bmod (\sqrt{-3})^4$, i.e. mod 9.
    \label{Cubes in quad extension} 
\end{lemma}
\begin{proof}
    Let $c \in \Q_3(\sqrt{-3})^\times$ be a cube. The first part is obvious, i.e., $v_{\sqrt{-3}}(c) = 3m$. Now without loss of generality let $c \in \Z_3[\sqrt{-3}]^\times$. Then for some $a,b \in \Z_3$ (where $a$ nonzero)
    \begin{align*}
        c &= (a + b\sqrt{-3})^3 \\
        &= a^3 + 3a^2b\sqrt{-3} - 9ab^2 - 3b^3\sqrt{-3} \\
        &= a^3 + 3(a-b)(a+b)b\sqrt{-3} - 9ab^2
    \end{align*}
    We have $v_3(a) = 0$, so $a \equiv 1$ or $-1 \bmod 3$. It's not hard to see from this that one of $b, a +b, a-b$ is $\equiv 0$ mod 3. So
    \[
        c \equiv a^3 \bmod (\sqrt{-3})^4
    \]
    But of course $(\sqrt{-3})^4 = 9$. So if $c \in \Z_3[\sqrt{-3}]^\times$ is a cube, it is equivalent mod 9 to the cube of some $a\in\Z_3^\times$. By Lemma \ref{cubes in Q3}, this is equivalent to $c \equiv \pm 1 \bmod 9 = (\sqrt{-3})^4$.

    Conversely, if $c \equiv \pm 1 \bmod 9$, then the polynomial $x^3 - c$ has $\pm 1$ as a solution mod $(\sqrt{-3})^4$. Furthermore, the derivative at that solution is $3(\pm 1)^2 = 3 \not\equiv 0 \bmod (\sqrt{-3})^4$. So by Hensel's lemma this lifts to a solution of $x^3 - c$ in $\Z_3(\sqrt{-3})$, and $c$ is a cube.
\end{proof}
\begin{corollary}\label{cor:quotient-norm-group-triviality-criterion}
    $c \in \Q_3(\sqrt{-3})^\times$ is in $\Q_3^\times(\Q_3(\sqrt{-3})^\times)^3$ iff its unit part $u$ is equivalent mod $(\sqrt{-3})^4$ to an invertible integer (i.e. $\{1,2,4,5,7,8\}$).
\end{corollary}
\begin{proof}
    First note that $\sqrt{-3} = -\frac{1}{3}\cdot (\sqrt{-3})^3 \in \Q_3^\times(\Q_3(\sqrt{-3})^\times)^3$, so $c \in \Q_3^\times(\Q_3(\sqrt{-3})^\times)^3$ iff $u$ is, and we can work directly with $u$.

    Then the if direction is obvious. For the other direction, suppose there exists $d \in \Q_3^\times$ such that $du \in (\Q_3(\sqrt{-3})^\times)^3$. Then $du = ({\sqrt{-3}})^{3m} w$ where $w\in \Z_3(\sqrt{-3})^\times$, $w \equiv \pm 1 \bmod (\sqrt{-3})^4$. Then $u \equiv d\inv ({\sqrt{-3}})^{3m} w \bmod (\sqrt{-3})^4$, where $v_{\sqrt{-3}}(d\inv ({\sqrt{-3}})^{3m}) = 0$, so $m$ even and $d\inv ({\sqrt{-3}})^{3m} \in \Z_3^\times$. Thus mod 9, $d\inv ({\sqrt{-3}})^{3m}$ is an invertible integer, and so is $u$.
\end{proof}

In the following proposition, we produce the nontrivial elements of the reduced cubic norm group for each non-Galois case. 
\begin{proposition} \label{prop:computing reduced cubic norm groups in non-galois}
For $\tau \in \{1,4,7\}$, let $\theta_\tau$ a root of $f_\tau(x) = x^3 + 3\tau$. Then the reduced cubic norm group of $\Q_3(\theta_\tau, \sqrt{-3})/\Q_3(\sqrt{-3})$ is represented (non-redundantly) in $\Q_3(\sqrt{-3})^\times/\Q_3^\times(\Q_3(\sqrt{-3})^\times)^3$ 
by the elements
\[
    \left\{1, \tau + \sqrt{-3}, \tau^2 - 3 + 2\tau\sqrt{-3}\right\}
\]

Let $\gamma$ a root of $g(x) = x^3 + 3x + 3$. Then the reduced cubic norm group of $\Q_3(\gamma, \sqrt{-3})/\Q_3(\sqrt{-3})$ is represented non-redundantly by
\[
\left\{1, 1 + (\sqrt{-3})^3, 1 + 2(\sqrt{-3})^3\right\}
\]

\end{proposition}
\begin{proof}
    We know from Lemma \ref{reduced cubic norm group isom to C3} that the reduced cubic norm group is $\cong C_3$.
 
    In the $x^3 + 3\tau$ cases, we compute $N(\theta_\tau + \sqrt{-3}) = -3\tau + (\sqrt{-3})^3$. Since $v_{\sqrt{-3}}(\tau) = 0$, the unit part of this is $\tau + \sqrt{-3}$. This is non-integer mod 9, so by Corollary \ref{cor:quotient-norm-group-triviality-criterion} it is not in $\Q_3^\times (\Q_3(\sqrt{-3})^\times)^3$, and thus represents a nontrivial element of the reduced cubic norm group. Squaring it and noticing that 

    In the $x^3+3x+3$ case,
    \begin{equation*}
        N(\gamma^2\sqrt{-3} + \gamma + 1) = 1 + (\sqrt{-3})^3 + (\sqrt{-3})^7 \equiv 1 + (\sqrt{-3})^3 \mod (\sqrt{-3})^4
    \end{equation*}
    so $\frac{N(\gamma^2\sqrt{-3} + \gamma + 1)}{1 + (\sqrt{-3})^3} \equiv 1 \bmod 9$. Then by Corollary \ref{cor:quotient-norm-group-triviality-criterion} they represent the same element of $\Q_3(\sqrt{-3})^\times/\Q_3^\times(\Q_3(\sqrt{-3})^\times)^3$, and since $1 + (\sqrt{-3})^3$ is not $\equiv$ an integer mod $(\sqrt{-3})^4$, this element is a nontrivial one. Then
    \[
    (1 + (\sqrt{-3})^3)^2 = 1 + 2(\sqrt{-3})^3 + (\sqrt{-3})^6 \equiv 1 + 2(\sqrt{-3})^3 \mod (\sqrt{-3})^4
    \]
    so $(1 + (\sqrt{-3})^3)^2$ and $1 + 2(\sqrt{-3})^3$ both represent the other nontrivial element of the reduced cubic norm group.
\end{proof}

\subsection{Determining the isomorphism class}

Now we are ready to compare the norm group of the given polynomial and the 4 non-Galois, $\sqrt{3c}$ extensions.

\begin{proposition}\label{prop:non-Galois reduced norm group generator}
    Let $f(x) = x^3 + \alpha x + \beta$ a depressed irreducible cubic with $\frac{\Delta_f}{-3}$ square in $\Q_3^\times$, $\theta$ a root of $f$.
    
    Let $u$ the unit part of $\beta$, $w$ the unit part of $\alpha$. Let $r = v_3(\alpha)$, $m = \lfloor v_3(\beta)/3 \rfloor$. Then
    \begin{align*}
        \begin{cases}
            1 + u(\sqrt{-3})^3 & v_3(\beta) \equiv 0 \bmod 3\\
            1 - 3^{r - 2m - 1} w + u\sqrt{-3} & v_3(\beta) \equiv 2 \bmod 3\\
            u + \left( 1 - 3^{r - 2m - 1} w\right) \sqrt{-3} & v_3(\beta)\equiv 1 \bmod 3,\quad r > 2m + 1\\
            1 + 2w (\sqrt{-3})^3 & v_3(\beta)\equiv 1 \bmod 3,\quad r = 2m + 1
        \end{cases}
    \end{align*}
    represents a nontrivial element of the reduced cubic norm group of $\Q_3(\theta, \sqrt{-3})/\Q_3(\sqrt{-3})$. In the second case, $r-2m \geq 1$, and in the third $r-2m\geq 2$.
\end{proposition}
\begin{proof}
    Recall that by Corollary \ref{cor:quotient-norm-group-triviality-criterion}, an element of $\Q_3(\sqrt{-3})$ is in $\Q_3^\times(\Q_3(\sqrt{-3})^\times)^3$ iff its unit part is equivalent mod 9 to an invertible integer (i.e. $\{1,2,4,5,7,8\}$).

    Also recall that $-1$, $\sqrt{-3} \in \Q_3^\times (\Q_3(\sqrt{-3})^\times)^3$, so we can multiply by them freely.

    We have by assumption that
    \begin{align*}
        \beta &= 3^{3m+\ell} u = (-1)^{m+\ell} u(\sqrt{-3})^{6m + 2\ell}\\
        \alpha &= 3^r w = (-1)^r w (\sqrt{-3})^{2r}
    \end{align*}
    By Lemma \ref{lem:wild-newton-polygon}, $v_3(\alpha) - \frac{2}{3} v_3(\beta) > 0$. We can write $v_3(\beta) = 3m + \ell$, where $\ell \in \{0,1,2\}$. Then $r - 2m >  2\ell/3$.

    Now, for any $\lambda \in \Q_3$, we can calculate
    \begin{align*}
        N(\theta + \lambda) &= -\beta + \lambda^3 + \lambda\alpha\\
        N(\theta^2 + \lambda) &= -\beta^2 + \lambda^3 - 2\lambda^2\alpha + \lambda\alpha^2
    \end{align*}
    We'll use these in each case to find norms that are nontrivial in the reduced cubic norm group.

    \bigskip

    Suppose $\ell = 0$. Then $r - 2m >  0$, so $r-2m \geq 1$.  Take $\lambda = (-1)^{m+1}(\sqrt{-3})^{2m-1}$. Then
    \begin{align*}
        (-1)^{m+1}(\sqrt{-3})^{-6m+3} N(\theta + \lambda)
        &= (-1)^{m+1}(\sqrt{-3})^{-6m+3} (-\beta + \lambda^3 + \lambda\alpha)\\
        &= (-1)^{m+1}\left(-(-1)^m u(\sqrt{-3})^{3} + (-1)^{3m +3} y^3 + (-1)^{m+1+r} (\sqrt{-3})^{2r + 2m - 1 - 6m + 3}\right)\\
        &= u(\sqrt{-3})^{3} + 1 +  (-1)^{r} w(\sqrt{-3})^{2(r - 2m) + 2}
    \end{align*}
    where $2(r - 2m) + 2 \geq 4$. So this is a unit, and $\equiv 1 + u(\sqrt{-3})^3 \bmod (\sqrt{-3})^4$, which is not an integer mod $(\sqrt{-3})^4$, since it has a term with $v_{\sqrt{-3}} = 3$. Thus by \ref{cor:quotient-norm-group-triviality-criterion} it is nontrivial.

    \bigskip

    Suppose $\ell = 2$. Then $r - 2m >  4/3$, so $r-2m \geq 2$. Take $\lambda = (-1)^{m+1} (\sqrt{-3})^{2m + 1}$. Then
    \begin{align*}
        (-1)^{m+1}(\sqrt{-3})^{-6m-3} N(\theta + \lambda)
        &= (-1)^{m+1} (\sqrt{-3})^{-6m-3} (-\beta + \lambda^3 + \lambda\alpha)\\
        &= (-1)^{m+1} (-(-1)^{m+2} u \sqrt{-3} + (-1)^{3m+3} + (-1)^{m+1+r} w (\sqrt{-3})^{2r+2m + 1 - 6m - 3})\\
        &= u\sqrt{-3} + 1 + (-1)^r w(\sqrt{-3})^{2(r - 2m) - 2}\\
        &= 1 - 3^{r - 2m - 1} w + u\sqrt{-3}
    \end{align*}
    where $r - 2m - 1 \geq 1$. So this is a unit, and since it has a nonzero term with $v_{\sqrt{-3}} = 1$, it is noninteger mod $(\sqrt{-3})^4$ and by \ref{cor:quotient-norm-group-triviality-criterion} nontrivial.

    \bigskip
    
    Suppose $\ell = 1$. Then $r - 2m >  2/3$, so $r-2m \geq 1$.
    
    If $r - 2m > 1$, take $\lambda = (-1)^m (\sqrt{-3})^{2m + 1}$. Then
    \begin{align*}
        (-1)^m (\sqrt{-3})^{-6m - 2} N(\theta + \lambda)
        &= (-1)^m (\sqrt{-3})^{-6m - 2} (-\beta + \lambda^3 + \lambda\alpha)\\
        &= (-1)^m(-(-1)^{m+1} u + (-1)^{3m} \sqrt{-3} + (-1)^{m+r} w (\sqrt{-3})^{2r + 2m + 1 - 6m - 2}\\
        &= u + \sqrt{-3} + (-1)^r w (\sqrt{-3})^{2(r - 2m) - 1}\\
        &= u + \left( 1 - 3^{r - 2m - 1} w\right) \sqrt{-3}
    \end{align*}
    If $r-2m>2$, this is $\equiv u + \sqrt{-3}$ mod $(\sqrt{-3})^4$. If $r-2m=2$, it is $\equiv u + \sqrt{-3} + w(\sqrt{-3})^3$. Either way, it has a nonzero term with $v_{\sqrt{-3}} = 1$, so is a noninteger and thus by \ref{cor:quotient-norm-group-triviality-criterion} nontrivial. 

    \medskip
    
    Finally, let $r-2m=1$, and take $\lambda = (\sqrt{-3})^{4m - 1}$. Then
    \begin{align*}
        (\sqrt{-3})^{-12m + 3} N(\theta^2 + \lambda)
        &= (\sqrt{-3})^{-12m + 3} (-\beta^2 + \lambda^3 - 2\lambda^2\alpha + \lambda\alpha^2)\\
        &= (-u^2 (\sqrt{-3})^{7} + 1 + 2(-1)^{r+1} w (\sqrt{-3})^{2r + 8m - 2 - 12m + 3} + w^2 (\sqrt{-3})^{4r + 4m - 1 - 12m + 3})\\
        &= (-u^2 (\sqrt{-3})^{7} + 1 + 2w (\sqrt{-3})^{3} + w^2 (\sqrt{-3})^{6})
    \end{align*}
    which is $\equiv 1 + 2w (\sqrt{-3})^3 \bmod (\sqrt{-3})^4$, which has a term with $v_{\sqrt{-3}} = 3$ and is thus nontrivial.
\end{proof}

Then we distinguish between these four cases:

\begin{theorem}
\label{main theorem non-Galois}
    Let $f(x) = x^3 + \alpha x + \beta$ be a depressed irreducible cubic such that $\frac{\Delta_f}{-3}$ is a square in $\Q_3$. Let $\theta$ a root of $f$.

    Let $u$ the unit part of $\beta$, $w$ the unit part of $\alpha$ (if $\alpha\neq 0$), $r = v_3(\alpha)$, $m = \lfloor v_3(\beta)/3 \rfloor$. Then
    \begin{enumerate}
        \item If $v_3(\beta) \equiv 0 \bmod 3$, $\Q_3(\theta)/\Q_3 \cong \Q_3[x]/(x^3 + 3x + 3)$.

        \item If $v_3(\beta) \equiv 2 \bmod 3$, $\Q_3(\theta)/\Q_3 \cong \Q_3[x]/(x^3 + 3\tau)$ for the unique $\tau \in \{1,4,7\}$ such that
        \[
            \pm\tau \equiv \frac{1 + w(-3)^{r-2m-1}}{u} \equiv \begin{cases}
                \frac{1-3w}{u} & r = 2m + 2\\
                \frac{1}{u} & r > 2m + 2
            \end{cases} \quad\mod 9
        \]

        \item If $v_3(\beta) \equiv 1 \bmod 3$ and $r > 2m + 1$, $\Q_3(\theta)/\Q_3 \cong \Q_3[x]/(x^3 + 3\tau)$ for the unique $\tau \in \{1,4,7\}$ such that
        \[
            \pm\tau \equiv \frac{u}{1 - 3^{r - 2m - 1} w} \equiv \begin{cases}
                \frac{u}{1-3w} & r = 2m + 2\\
                u & r > 2m + 2
            \end{cases} \quad \mod 9
        \]

        \item If $v_3(\beta) \equiv 1 \bmod 3$ and $r = 2m + 1$, $\Q_3(\theta)/\Q_3 \cong \Q_3[x]/(x^3 + 3x + 3)$.
    \end{enumerate}
    
\end{theorem}
\begin{proof}
    By Lemma \ref{lem:non-galois discriminant criterion}, we know that $\Q_3(\theta)/\Q_3$ is isomorphic to the extension generated by $x^3 + 3\tau$ for $\tau \in \{1,4,7\}$, or by $x^3 + 3x + 3$.
    
    By Proposition \ref{prop:computing reduced cubic norm groups in non-galois}, in the first case its reduced cubic norm group is represented by $\{1, \tau + \sqrt{-3}, \tau^2 - 3 + 2\tau \sqrt{-3}\}$, and in the second case by $\{1, 1 + (\sqrt{-3})^3, 1 + 2(\sqrt{-3})^3\}$.
    
    \begin{enumerate}
        \item If $v_3(\beta) \equiv 0 \bmod 3$, then by Proposition \ref{prop:non-Galois reduced norm group generator}, $1 + u(\sqrt{-3})^3$ represents a nontrivial element of the reduced cubic norm group.

        Then 
        \[
        1 + u(\sqrt{-3})^3 \equiv \{1 + (\sqrt{-3})^3, 1 + 2(\sqrt{-3})^3\} \mod 9
        \]
        depending on whether $u \equiv 1$ or $2$ mod 3. So the reduced cubic norm group associated with $f$ is the same as that associated with $x^3 + 3x + 3$, and thus they generate isomorphic extensions.

        \item If $v_3(\beta) \equiv 2 \bmod 3$, Proposition \ref{prop:non-Galois reduced norm group generator} gives $1 + w(-3)^{r-2m-1} - u\sqrt{-3}$, where $r-2m \geq 2$.

        We can immediately rule out $x^3 + 3x + 3$, since the representatives of its reduced cubic norm group don't contain any nonzero term with $v_{\sqrt{-3}} = 1$ (i.e. a term of the form $y\sqrt{-3}$ where $y$ a unit).

        To distinguish between the $x^3 + 3\tau$ cases, note that the ratio mod 9 between the coefficients of the $v_{\sqrt{-3}} = 0$ and $v_{\sqrt{-3}} = 1$ terms is $\tau$ for $\tau + \sqrt{-3}$ and
        \[
        \frac{\tau^2 - 3}{2\tau} \equiv -\tau \bmod 9
        \]
        for $\tau^2 - 3 + 2\tau \sqrt{-3}$.

        Furthermore, this ratio is preserved when we multiply by a unit in $\Q_3^\times (\Q_3(\sqrt{-3})^\times)^3$, since by Corollary \ref{cor:quotient-norm-group-triviality-criterion} those are just invertible integers mod 9. So all unit representatives of nontrivial elements of the reduced cubic norm group share this ratio, up to sign.

        In particular, since the extension generated by $f$ is isomorphic to one generated by $x^3 + 3\tau$ for some (unique) $\tau$, this means that the ratio between the integer and $\sqrt{-3}$ terms of $1 + w(-3)^{r-2m-1} - u\sqrt{-3}$ is $\pm \tau$ for that $\tau$ exactly.

        Then we can compute this ratio, and thus $\pm \tau$, as (ignoring sign)
        \[
        \frac{1 + w(-3)^{r-2m-1}}{u} \equiv \begin{cases}
        \frac{1-3w}{u} & r = 2m + 2\\
        \frac{1}{u} & r > 2m + 2
        \end{cases} \quad\mod 9
        \]

        \item If $v_3(\beta) \equiv 1 \bmod 3$ and $r > 2m + 1$, Proposition \ref{prop:non-Galois reduced norm group generator} gives $u + \left( 1 - 3^{r - 2m - 1} w\right) \sqrt{-3}$.

        Then by the same argument as in the previous case, we can rule out $x^3 + 3x + 3$ and determine  $\pm\tau$ mod 9 by taking the ratio between the two coefficients, which is
        \[
        \frac{u}{1 - 3^{r - 2m - 1} w} \equiv \begin{cases}
            \frac{u}{1-3w} & r = 2m + 2\\
            u & r > 2m + 2
        \end{cases} \quad \mod 9
        \]

        \item If $v_3(\beta) \equiv 1 \bmod 3$ and $r = 2m + 1$, Proposition \ref{prop:non-Galois reduced norm group generator} gives $1 + 2w (\sqrt{-3})^3$. As in the first case, we see immediately that
        \[
        1 + 2w (\sqrt{-3})^3 \equiv \{1 + (\sqrt{-3})^3, 1 + 2(\sqrt{-3})^3\} \mod 9
        \]
        and therefore the extension is isomorphic to that generated by $x^3 + 3x + 3$.

    \end{enumerate}
\end{proof}

\section{Acknowledgements}

The authors would like to acknowledge the invaluable guidance provided by our project supervisors, Daniel Johnstone, Camelia Karimianpour, and Malors Espinosa Lara. We are especially grateful to Daniel Johnstone for his key insights and detailed feedback on early drafts of this paper. This project was generously supported by the Fields Institute for Research in Mathematical Sciences through the 2023 FUSRP summer research program.

\bibliographystyle{amsalpha}
\bibliography{bibliography}

\end{document}